\documentclass[12pt]{amsart}
\usepackage{amscd,amssymb,amsmath}
\begin{document}

\def\supp{\operatorname{supp}}
\def\Ind{\operatorname{Ind}}
\def\interior{\operatorname{int}}

\def\Aut{\operatorname{Aut}}
\def\Ad{\operatorname{Ad}}
\def\range{\operatorname{range}}
\def\sp{\operatorname{sp}}
\def\clsp{\overline{\operatorname{sp}}}
\def\Prim{\operatorname{Prim}}
\def\dashind{\operatorname{\!-Ind}}
\def\id{\textnormal{id}}
\def\rt{\textnormal{rt}}
\def\lt{\textnormal{lt}}

\def\H{\mathcal{H}}
\def\L{\mathcal{L}}
\def\K{\mathcal{K}}

\def\B{\mathbf{B}}
\def\C{\mathbb{C}}
\def\R{\mathbb{R}}
\def\Z{\mathbb{Z}}
\def\T{\mathbb{T}}

\hyphenation{iso-metry}

\newtheorem{thm}{Theorem}[section]

\newtheorem{cor}[thm]{Corollary}
\newtheorem{lem}[thm]{Lemma}
\newtheorem{lemma}[thm]{Lemma}
\newtheorem{prop}[thm]{Proposition}

\theoremstyle{definition}
\newtheorem{definition}[thm]{Definition}

\theoremstyle{remark}
\newtheorem{remark}[thm]{Remark}
\newtheorem{ex}[thm]{Example}
\newtheorem{example}[thm]{Example}
\newtheorem{remarks}[thm]{Remarks}
\newtheorem{claim}[thm]{Claim}
\newtheorem{problem}[thm]{Problem}

\numberwithin{equation}{section}

\title[Proper actions]{\boldmath Proper actions on imprimitivity bimodules\\
and decompositions of Morita equivalences}

\author[an Huef]{Astrid an Huef}
\address{Department of Mathematics and Computer Science\\
University of Denver\\
2360 S.~Gaylord St.\\
Denver, CO 80208-0189\\
USA}
\email{astrid@cs.du.edu}

\author[Raeburn]{Iain Raeburn}
\address{Department of Mathematics\\ University of Newcastle\\
NSW 2308\\ Australia}
\email{iain@frey.newcastle.edu.au}

\author[Willimas]{Dana P. Williams}
\address{Department of Mathematics\\
  Dartmouth College\\
  Hanover, NH 03755-3551\\
  USA} \email{dana.williams@dartmouth.edu}

\thanks{This research was supported by grants from the University of Denver, the University of Newcastle, and Dartmouth College. }

\date{April 10, 2001}

\begin{abstract}
We consider a class of proper actions of locally compact groups on imprimitivity bimodules over $C^*$-algebras which behave like the proper actions on $C^*$-algebras introduced by Rieffel in 1988. We prove that every such action gives rise to a Morita equivalence between a crossed product and a generalized fixed-point algebra, and in doing so make several innovations which improve the applicability of Rieffel's theory. We then show how our construction can be used to obtain canonical tensor-product decompositions of important Morita equivalences. Our results show, for example, that the different proofs of the symmetric imprimitivity theorem for actions on induced algebras yield isomorphic equivalences. A similar analysis of the symmetric imprimitivity theorem for graph algebras gives new information about the amenability of actions on graph algebras. 
\end{abstract}

\maketitle

\section{Introduction}

A Morita equivalence between two $C^*$-algebras $A$ and $B$ is
implemented by an imprimitivity bimodule ${}_AX_B$, which carries the
structure necessary to induce Hilbert space representations from $B$ to
$A$ and back again. There are often several ways of constructing these
bimodules, and, unsurprisingly, some ways are better for some things, and
others for others. One therefore wants to be able to switch between
different bimodules implementing equivalences between the same algebras.

To illustrate the kind of problems which arise, we consider a situation
which underlies many important equivalences. Suppose we have commuting free and
proper actions of locally compact groups $K$ and $H$ on the left and
right of a locally compact Hausdorff space $P$. The orbit spaces
$P/H$ and $K\backslash P$  are again locally compact Hausdorff spaces, and
carry actions of, respectively, $K$ and $H$; the symmetric imprimitivity theorem of Green
and Rieffel says that the crossed products $C_0(P/H)\rtimes K$ and
$C_0(K\backslash P)\rtimes H$ are Morita equivalent. In the original
proof of \cite{rieff-applications}, a suitable imprimitivity bimodule $W$
was constructed by completing the space $C_c(P)$ of continuous functions of compact 
support. It was later shown in \cite{cmw}
that one could appeal to a previous construction of Green which gives
$C_0(P/H)$--$(C_0(P)\rtimes H)$ and $(C_0(P)\rtimes K)$--$C_0(K\backslash
P)$ bimodules $X$ and $Y$, form crossed product bimodules $X\rtimes K$
and $Y\rtimes H$, and take the internal tensor product $(X\rtimes
K)\otimes_{C_0(P)\rtimes(K\times H)}(Y\rtimes H)$ as the desired
$(C_0(P/H)\rtimes K)$--$(C_0(K\backslash P)\rtimes H)$ imprimitivity
bimodule. This latter construction has advantages: for example, it saves
burrowing into the detailed construction of bimodules, allows us to
analyze the effect of extra structure in stages, and makes it easier to
prove analogues for reduced crossed products. On the other hand, we have
available a concrete bimodule $W$, which is  much more convenient for
direct calculations. To make the best of both worlds, we need to prove
that
\begin{equation}\label{fundiso}
W\cong (X\rtimes
K)\otimes_{C_0(P)\rtimes(K\rtimes H)}(Y\rtimes H) 
\end{equation}
as $(C_0(P/H)\rtimes K)$--$(C_0(K\backslash P)\rtimes H)$ imprimitivity
bimodules.

We ran into  problems like these in \cite{hrw}, where we found an isomorphism implementing~\eqref{fundiso} using {\it ad hoc} methods; to verify that it worked, we had to
do some awful calculations involving quintuple integrals. One goal of the present project was
to find a more systematic way of identifying and verifying such
isomorphisms: our  Theorem~\ref{thm-decomposition} tells us
not just that there is an isomorphism, but also how to write it down. 

To make our approach as systematic as possible, we have worked within the general framework of proper actions of groups on $C^*$-algebras, as developed by Rieffel in \cite{rie-pr}, and we have, we hope, made significant improvements to that theory. In particular, we have extended Rieffel's main Morita equivalence  in \cite[Corollary~1.7]{rie-pr} to cover proper actions on imprimitivity bimodules. This extension turned out to be anything but routine, and we are optimistic that some of the technical tools we have developed will help in constructing Morita equivalences  for more general integrable actions, where substantial technical problems arise (see \cite[\S6]{rie-pr2}).
Because Rieffel's framework involves reduced crossed products rather than full ones, our main results are about reduced crossed products. We intend to apply our techniques to full crossed products elsewhere.

After a brief first section in which we review the necessary integration theory, we begin in \S\ref{genRieffel} by discussing proper actions on imprimitivity bimodules. We start with a Morita
equivalence $(X,G,\gamma)$ between two $C^*$-dynamical systems $(A,G,\alpha)$ and  $(B,G,\beta)$. The action $\gamma$ is \emph{proper} if there is a $\gamma$-invariant pre-imprimitivity bimodule ${}_{A_0}(X_0)_{B_0}\subset X$  with properties like those of the dense subalgebra used by Rieffel. There are several ways in which the technical hypotheses could be phrased; we have chosen one which reduces to that of \cite{rie-pr} 
when ${}_AX_B={}_BB_B$, 
bears a striking formal resemblance to it, and yields the desired Morita equivalence $\overline{X_0}$ between $A\rtimes_{\alpha,r}G$ and a generalized fixed-point algebra $B^\beta$ when the action is also saturated (Theorem~\ref{thm-sumup}). The proof of Theorem~\ref{thm-sumup}, though, is quite different from its analogue in \cite{rie-pr}. For ${}_AX_B={}_BB_B$, Rieffel proved that the $(B\rtimes_{\beta,r} G)$-valued
inner product has the required properties, that $B^\beta$ acts as
adjointable operators on the resulting left Hilbert module ${}_{B\rtimes_r G}Z$, and then
that the map $B^\beta\to\L({}_{B\rtimes_r G}Z)$ is isometric \cite[page~151]{rie-pr}. We were not able to extend this last part, so we had to
substantially reshape the whole argument, starting with the right inner
product rather than the left. In retrospect, this is probably a good
thing. The process we have gone through is similar to the program
discussed by Rieffel in his later paper \cite[\S6]{rie-pr2}, and since
we have been able to sidestep some of the general problems he raises in
our setting, our arguments may be useful in the more general context. Indeed, we have already used
some of these ideas to find new insight on how the symmetric imprimitivity theorem relates to reduced crossed products (see \cite{aHRqs}). 

In \S\ref{decomp}, we prove our general decomposition theorem. The key idea is that one obtains a decomposition like (\ref{fundiso}) whenever one has a Morita equivalence for the linking algebra $L(X)$ of another Morita equivalence; the key Lemma~\ref{lem-ibm-isom} is a one-sided version of a result from \cite{es}. The main work in \S\ref{decomp} is to show that if the action $\gamma$  of $G$ on $X$ is proper and saturated, then so are the associated actions on $B$ and $X\oplus B$; we then apply Lemma~\ref{lem-ibm-isom} to a bimodule over $L(X\rtimes G)$ arising from an application of Theorem~\ref{thm-sumup} to $X\oplus B$. The result is a tensor-product decomposition for the bimodule $\overline{X_0}$ of Theorem~\ref{thm-sumup}, which in the situation of the symmetric imprimitivity theorem turns out to be the desired isomorphism~(\ref{fundiso}).

We discuss the application to the symmetric imprimitivity theorem in \S\ref{sec-sit}. We work in the generality of \cite{rae}, considering actions of two groups $K$ and $H$ on different induced $C^*$-algebras. These induced algebras are the generalized fixed-point algebras for the diagonal actions considered in \cite{rw85}, and we apply Theorem~\ref{thm-decomposition} to the bimodule constructed in \cite{rw85}. We then identify the two components in the resulting decomposition as crossed products of two such bimodules. While this process is by no means trivial,   our theory  tells us where to go at each stage, and hence much of \S\ref{sec-sit} is therefore routine. In the last section, we give another application of Theorem~\ref{thm-decomposition} to crossed products of graph algebras. This gives new information about the symmetric imprimitivity theorem for graph algebras, and allows us to settle a question left open in  \cite{PR}.

\section{Background on integration}\label{intro-waffle}

Let $G$ be a locally compact group and $A$ a $C^*$-algebra.
We need to know that a continuous
function $f:G\to A$ satisfying $\int_G\|f(s)\|\,ds<\infty$ has a well-defined
integral  in $A$, and of course there is no problem with
this. For example, minor modifications to the construction of
\cite[Lemma~C.3]{tfb} give a satisfactory integral which is characterized
by its behaviour  under representations: whenever $\pi:A\to B(\H)$ is a
nondegenerate representation and $h,k\in \H$, the function $s\mapsto
\big(\pi(f(s))h\,\big|\,k\big)$ is integrable in the usual sense, and we
have
\[
\bigg(\pi\Big(\int_G f(s)\,ds\Big) h\;\bigg|\;k\bigg)
=\int_G\big(\pi(f(s))h\,\big|\,k\big)\,ds.
\]
It then follows that bounded linear maps and multiplication by
multipliers pull in and out of integrals with no problems. We also need to integrate continuous functions with values in Hilbert
modules, and the next lemma says we can do this in a way consistent with
our
$C^*$-algebra-valued integrals.

\begin{lemma}
Suppose $X$ is a right Hilbert $B$-module, and $f:G\to X$ is a continuous
map such that $\int_G\|f(s)\|_B\,ds<\infty$. Then there is a unique element
$\int_G f(s)\,ds$ of $X$ such that
\[
\big\langle x\,,\,{\textstyle \int_G f(s)\,ds}\big\rangle_B=
\int_G\langle x\,,\,f(s)\rangle_B\,ds.
\]
\end{lemma}

\begin{proof}
By viewing $X$ as the top right-hand corner in the linking algebra
$L(_{\K(X)}X_B)$, we can convert $f$ into a function with
values in a $C^*$-algebra, and integrate as usual. The inner product on
$X$ is given by multiplication in $L(X)$, and hence can be pulled in and
out of the integral.
\end{proof}

In \S\ref{sec-sit}, we shall want to pull variables through integrals:
for example, it is useful to know that the product $f*g$ of two functions in 
$C_c(G,A)\subset A\rtimes_\alpha G$ is given by the $(A\rtimes_\alpha G)$-valued 
integral $\int_G f(r)(\lt\otimes\alpha)_r(g)\,dr$ as well as the function 
$s\mapsto \int_G f(r)\alpha_r(g(r^{-1}s))\,ds$. The following lemma makes this 
possible.

\begin{lemma}\label{pullpin}
Suppose $(A,H,\alpha)$ is a dynamical system, $G$ is another locally compact 
group, and $f:G\times H\to A$ is a continuous function of compact support.
Then the map $F:r\mapsto f(r,\cdot)$ is a continuous function of compact 
support with values in $A\rtimes_\alpha H$, and $\int_G F(r)\,dr$ is represented
by the function $s\mapsto \int_G f(r,s)\,dr$ in $C_c(G,A)$. In other words,
\[
\Big(\int_G F(r)\,dr\Big)(s)=\int_G f(r,s)\,dr.
\]
\end{lemma}

\begin{proof}
Let $L$ be a compact neighbourhood of $\{s:f(r,s)\not=0\mbox{ for some }r\}$, 
and denote its interior by $\interior L$. Then by uniform continuity, the function $F$ 
is continuous from $G$ into the $C^*$-algebra $C_0(\interior L,A)$. Thus it has an
integral
\[
\int_G^{C_0(\interior L, A)} F(r)\,dr
\] 
in $C_0(\interior L,A)$; since evaluation at $s$ is a continuous homomorphism on
$C_0(\interior L,A)$, we have
\begin{equation}\label{evalonC0}
\Big(\int_G^{C_0(\interior L, A)} F(r)\,dr\Big)(s)=\int_G^A f(r,s)\,dr
\end{equation}
for every $s\in G$. The inclusion $i$ of $C_0(\interior L,A)\subset C_c(G,A)$ into
$A\rtimes_\alpha G$ is a bounded linear map, and hence
\begin{equation}\label{pullithru}
\int_G^{A\rtimes_\alpha G} F(r)\,dr= \int_G^{A\rtimes_\alpha G} i(F(r))\,dr=
i\Big(\int_G^{C_0(\interior L, A)} F(r)\,dr\Big).
\end{equation}
Putting (\ref{evalonC0}) and (\ref{pullithru}) together says that the integral
$\int_G^{A\rtimes_\alpha G} F(r)\,dr$ is given by the continuous function of 
compact support defined by the right-hand side of~\eqref{evalonC0}.
\end{proof}

\begin{remark}\label{pullpinbimods}
We will also want to pull the variable through integrals taking values in 
imprimitivity bimodules $X$ obtained by completing $C_c(P,A)$, where $P$ 
is a locally compact space and $A$ is a $C^*$-algebra. The above argument 
shows that this is valid provided the maps $i:C_0(\interior L,A)\to X$ are 
continuous for the sup-norm on $C_0(\interior L, A)$. This is the case 
for all the bimodules considered in \S\ref{sec-sit}.
\end{remark}

\section{Proper actions on imprimitivity bimodules}\label{genRieffel}

Throughout this section,
$(X,G,\gamma)$ will be a Morita equivalence between two dynamical systems
$(A,G,\alpha)$ and $(B,G,\beta)$; since there is only one locally
compact group $G$ involved, we will drop it from our notation. 

\begin{definition}
\label{defn-pr}
The action $\gamma$ of $G$ on ${}_A X_B$
is
\emph{proper} if there are an invariant subspace
$X_0$ of $X$ and invariant $*$-subalgebras $A_0$ of $A$
and $B_0$ of $B$, 
such that $_{A_0}(X_0)_{B_0}$ is a pre-imprimitivity bimodule with
completion $_AX_B$, and such that
\begin{enumerate}
\item[(1)]
for every $x,y\in X_0$, the functions 
$s\mapsto \Delta(s)^{-1/2}{}_A\langle
x\,,\,\gamma_s(y)\rangle$ and $s\mapsto 
{}_A\langle x\,,\,\gamma_s(y)\rangle$ are in
$L^1(G,A)$;
\item[(2)]
for every $b\in B_0$ and $x\in X_0$, the functions $s\mapsto \gamma_s(x)\cdot b$ and $s\mapsto \Delta(s)^{-1/2}\gamma_s(x)\cdot b$ are in
$L^1(G,X)$;
\item[(3)]for every $x,y\in X_0$, there is a multiplier
$\langle x\,,\,y\rangle_D$ in $M(B_0)^\beta$ such that $z\cdot\langle
x\,,\,y\rangle_D\in X_0$ for all $z\in X_0$, and
\begin{equation}\label{pr-eq}
\int_G b\beta_s(\langle x\,,\,y\rangle_B)\, ds
=b\cdot\langle x\,,\,y\rangle_D\ \mbox{ for all }b\in B_0.
\end{equation}
\end{enumerate}  
\end{definition}

That the integral in~\eqref{pr-eq} exists follows from 
\begin{equation*}\big\| b\beta_s\big(\langle x\,,\, y\rangle_B\big)\big\|=\big\|\langle\gamma_s(y)\,,\,\gamma_s(x)\cdot b^*\rangle_B\big\|
\leq \|\gamma_s(x)\cdot b^*\|\,\|y\|
\end{equation*}
and item (2) of the definition.

\begin{remarks} 
\begin{enumerate}
\item[(1)]
There are some subtleties to this definition.
First, asserting that
$_{A_0}(X_0)_{B_0}$ is a pre-imprimitivity bimodule is an efficient way
of saying many things; for example, it implies that $x\cdot b\in X_0$
whenever $x\in X_0$ and
$b\in B_0$, and
$_A\langle x\,,\,y\rangle\in A_0$ whenever $x,y\in X_0$. Second, saying
that $X$ is the completion of $X_0$ is meant to include that $A_0$
is dense in $A$ and $B_0$ is dense in $B$. And third, the $D$ adorning the
inner product does not yet exist: it will be defined in
Proposition~\ref{rtHstructure} below, and proved there that
$\langle\cdot,\cdot\rangle_D$ is a $D$-valued inner product.
\item[(2)] 
Note that the action $\beta$ of $G$ on ${}_BB_B$ is proper with respect to ${}_{B_0}(B_0)_{B_0}$ if and only if the action $\beta$ on $B$ is proper with respect to $B_0$ in the sense  \cite[Definition~1.2]{rie-pr}.
\item[(3)] 
Example~\ref{ex-green} shows that Definition~\ref{defn-pr} is not symmetric: asserting that $\gamma$ is proper is not the same as asserting that the action $\widetilde\gamma$ on the dual equivalence ${}_B\widetilde X_A$ is proper.
\end{enumerate}
\end{remarks}

We will prove that if the Morita equivalence ${}_{(A,\alpha)}(X,\gamma)_{(B,\beta)}$ is proper with respect to ${}_{A_0}(X_0)_{B_0}$ then $X_0$ completes to a Morita equivalence between an ideal $E$ of $A\rtimes_{\alpha,r}G$ and a generalized fixed-point algebra $D\subset M(B_0)^\beta$ of $B$. 

\begin{prop}\label{rtHstructure}
Let $D$ be the closure of $D_0:=\sp\{\langle x,y\rangle_D:x,y\in X_0\}$
in $M(B)$, where $\langle\cdot\,,\,\cdot\rangle_D$ is defined by
 Definition~\ref{defn-pr}(3). Then $D$ is a
$C^*$-algebra, and
$\langle\cdot\,,\,\cdot\rangle_D$ is a
$D$-valued inner product on $X_0$.
\end{prop}

\begin{proof}
Let $x,y\in X_0$ and $d\in D_0$. Then $d\in M(B_0)^\beta$, and so for every
$b\in B_0$, we have
\begin{align*}
b\langle x\,,\,y\rangle_D d
&=\Big(\int_G b\beta_t(\langle x\,,\,y \rangle_B)
\,dt \Big)d=\int_G b\beta_t(\langle x\,,\,y\rangle_B d)\,dt\\
&=\int_G b\beta_t(\langle x\,,\,y\cdot d\rangle_B)\,dt
= b\langle x\,,\,y\cdot d\rangle_D,
\end{align*} 
which implies $\langle x\,,\,y\rangle_Dd=\langle x\,,\,y\cdot
d\rangle_D$; since we know from Definition~\ref{defn-pr}(3) that $y\cdot
d\in X_0$, this implies that $D_0$ is an algebra.
Similarly, we have
\begin{align*}
b\langle y\,,\,x\rangle_D
&=\int_G b\beta_s(\langle y\,,\,x\rangle_B)\, ds
=\int_G b\beta_s(\langle x\,,\,y\rangle_B)^*\, ds\\
&=\Big(\int_G \beta_s(\langle x\,,\,y\rangle_B)b^* \, ds
\Big)^*
=(\langle x\,,\,y\rangle_D b^*)^*
=b\langle x\,,\, y\rangle_D^*,
\end{align*}
which implies that $\langle y\,,\,x\rangle_D=(\langle x\,,\,
y\rangle_D)^*$. This proves both that $D_0$ is a $*$-algebra, so its
closure $D$ is a $C^*$-algebra, and that
$\langle\cdot,\cdot\rangle_D$ has the algebraic properties of an inner
product.

To show positivity of $\langle\cdot,\cdot\rangle_D$, let $\pi$ be a
faithful nondegenerate  representation of $B$, and note that 
\begin{align*}
\langle x\,,\, x\rangle_D\geq 0 &\Longleftrightarrow\langle x\,,\, x\rangle_D\geq 0
\mbox{ in } M(B)\\
&\Longleftrightarrow(\bar\pi(\langle x\,,\,
x\rangle_D)h\,|\, h)\geq 0\mbox{ for all }h\in\H_\pi.
\end{align*} 
Since $B_0$ is dense in $B$ and $\pi$ is
nondegenerate, it is enough to show this when $h=\pi(b)k$ for $b\in
B_0$ and $k\in\H_\pi$. Well,
\begin{align*}
\big(\bar\pi(\langle x\,,\, x\rangle_D)\pi(b)k\,\big|\,\pi(b)k\big)
&=\big(\pi (b^* \langle x\,,\, x\rangle_D b)k\,\big|\,k\big)\\
&=\bigg(\pi\Big(\int_G b^*\beta_s(\langle 
x\,,\,x\rangle_B)b\,ds\Big)k\;\bigg|\; k\bigg)\\
&=\int_G \big(\pi(b^*\beta_s(\langle 
x\,,\,x\rangle_B)b)k\,\big|\, k\big)\,ds,
\end{align*}
which is positive because each  $b^*\beta_s(\langle
x\,,\,x\rangle_B)b$ is positive. Because $b^*\beta_s(\langle
x\,,\,x\rangle_B)b$ is continuous in $s$, this calculation also shows
that \[
\langle x\,,\, x\rangle_D=0\Longrightarrow\beta_s(\langle
x\,,\,x\rangle_B)=0 \mbox{ for all $s$ }\Longrightarrow x=0,
\]
so $\langle\cdot,\cdot \rangle_D$ is definite. 
\end{proof}

\begin{remark}
To prove that an element $d$ of a $C^*$-algebra $D$ is positive, it is
usually enough to take a faithful representation $\mu$ of $D$ and prove
that $\mu(d)\geq0$ as an operator on $\H_\mu$. For the above argument,
however, it is essential that the representation $\mu$ of $D_0$ is the
restriction of (an extension of) a representation $\pi$ of $B$. In general, not every
representation $\mu$ of $D_0$ arises this way, and choosing
$\mu=\bar\pi|_{D_0}$ means that we are proving $d\geq 0$ in the
$C^*$-algebra obtained by completing $D_0$ in the norm of
the $C^*$-algebra $M(B)$. This observation is crucial in \cite{aHRqs}.
\end{remark}

\begin{lemma}\label{leftaction}
There are homomorphisms $\mu:A\to \L(\overline{(X_0)_D})$ and $U:G\to
U\L(\overline{(X_0)_D})$ such that
\[
\mu(a)x=a\cdot x,\  U_sx=\Delta(s)^{1/2}\gamma_s(x)\ \mbox{and}\ 
\mu(\alpha_s(a))=U_s\mu(a)U_s^*\ \mbox{for $a\in A_0$, $x\in X_0$}.
\]
\end{lemma}

\begin{remark}\label{rem-cov?}
The homomorphism $\gamma:G\to\L(X_B)$ is not unitary: it
changes the inner product by $\beta_s$. So it is important here that we
are talking about the $D$-valued inner product on $X_0$, and that we
have had to introduce the modular function to ensure that $U_s$
preserves this inner product. At this stage, we are not asserting that
$(\mu,U)$ is a covariant representation in the usual sense: neither
nondegeneracy of
$\mu$ nor continuity of
$U$ seems obvious. We shall return to this
point in Lemma~\ref{formulaforaction} below.  Meanwhile, we observe that these two problems also arise in the construction of a Morita equivalence for more general integrable actions \cite[\S6]{rie-pr2}.
\end{remark}

\begin{proof}[Proof of Lemma~\ref{leftaction}]
We first show that $\langle a\cdot x,a\cdot x\rangle_D\leq \|a\|^2\langle
x\,,\,x\rangle_D$ as elements of $D$, so that $\mu(a):x\mapsto a\cdot x$
is bounded on $(X_0)_D$. To do this, we again choose a faithful
nondegenerate representation $\pi$ of $B$, and it is enough to prove that
\begin{equation}\label{mu(a)bded}
\big(\bar\pi\big(\|a\|^2\langle x\,,\,x\rangle_D- \langle a\cdot x,a\cdot
x\rangle_D\big)\pi(b)h\,\big|\,\pi(b)h\big)\geq 0
\end{equation}
for all $b\in B_0$ and $h\in \H_\pi$. We know that
\[
\|a\|^2\langle x\,,\,x\rangle_B- \langle a\cdot x\,,\,a\cdot
x\rangle_B
\]
is positive in $B$, so
\[
\big(\pi\big(b^*\beta_s(\|a\|^2\langle x\,,\,x\rangle_B- \langle a\cdot
x,a\cdot x\rangle_B)b\big)h\,\big|\,h\big)\geq 0
\]
for all $b$, $h$ and $s$. Integrating this over $G$ and pulling the
integral inside the inner product gives
\[
\bigg(\pi\Big(\int_G b^*\beta_s\big(\|a\|^2\langle x\,,\,x\rangle_B- \langle
a\cdot x\,,\,a\cdot x\rangle_B\big)b\,ds\Big)h\;\bigg|\;h\bigg)\geq 0,
\]
which is \eqref{mu(a)bded}. We deduce that $\mu(a)$
is bounded. Since $a^*$ is just another element of
$A_0$, 
$\mu(a^*)$ is also bounded.
We can therefore show that $\mu(a)$ is adjointable with adjoint
$\mu(a^*)$ by
checking that
\[
b\langle a\cdot x\,,\,y\rangle_D =b\langle x\,,\,a^*\cdot y\rangle_D\ \mbox{ for
$b\in B_0$},
\]
and this follows easily from $\langle a\cdot x,y\rangle_B =\langle
x,a^*\cdot y\rangle_B$. Thus $\mu(a)\in\L(\overline{(X_0)_D})$, and it is
easy to check that $\mu$ is a homomorphism of $C^*$-algebras.

To verify that $U_s$ is unitary, we let $b\in B_0$, $x,y\in X_0$ and
calculate:
\begin{align*}
b\langle U_sx\,,\,U_sy\rangle_D
&=\int_G b\beta_t\big(\langle
\Delta(s)^{1/2}\gamma_s(x)\,,\,\Delta(s)^{1/2}\gamma_s(y)\rangle_B\big)\,dt\\
&=\int_G b\beta_{ts}\big(\langle x\,,\,y \rangle_B\big)\Delta(s)\, dt\\
&=\int_G b\beta_{r}\big(\langle x\,,\,y \rangle_B\big)\,dr\\
&=b\langle x\,,\,y \rangle_D.
\end{align*}
This calculation shows that $U_s$ is bounded with $\|U_s\|=1$, as is $U_{s^{-1}}$, and that
\[
b\langle U_sx\,,\,y \rangle_D=b\langle U_sx\,,\,U_s(U_{s^{-1}}y) \rangle_D=
b\langle x\,,\,U_{s^{-1}}y \rangle_D,
\]
so $U_s$ is adjointable with $U_s^*=U_{s^{-1}}$. Since we trivially have
$U_sU_t=U_{st}$ on $X_0$, $U$ is a homomorphism into
$U\L(\overline{(X_0)_D})$, as claimed. 

For $x\in X_0$, the covariance condition
$\mu(\alpha_s(a))x=U_s\mu(a)U_s^*x$ follows easily from the formulas and
the identity $U_s^*=U_{s^{-1}}$, and it then extends by continuity to
$x\in \overline{(X_0)_D}$.
\end{proof}

As we observed in Remark~\ref{rem-cov?}, we do not know whether $(\mu,U)$ is always a covariant representation in
$\L(\overline{(X_0)_D})=M(\K(\overline{(X_0)_D}))$, but we can make it
into one by representing $\L(\overline{(X_0)_D})$ on Hilbert space. As
usual, we start with a nondegenerate faithful representation $\pi$ of
$B$. Then $\bar\pi$ is faithful on $D$, and hence
$X_0\dashind_D^\L\bar\pi$ is a faithful representation of
$\L(\overline{(X_0)_D})$ on $X_0\otimes_D\H_\pi$. We write $(\nu,V)$ for
the pair $\big((X_0\dashind_D^\L\bar\pi)\circ
\mu,(X_0\dashind_D^\L\bar\pi)\circ U\big)$, so that
\[
\nu(a)(x\otimes_D h)=(a\cdot x)\otimes_D h\ \mbox{ and }\ V_s(x\otimes_D
h)=\Delta(s)^{1/2}\gamma_s(x)\otimes_D h
\]
for $x\in X_0$ and $a\in A_0$. To see that $(\nu,V)$ is covariant, we
relate it to the right-regular representation
$\big((X\dashind_B^A\pi)\widetilde{\;},\rho\big)$ of
$(A,\alpha)$ on $L^2(G,X\otimes_B\H_\pi)$, which is given by
\[
(X\dashind_B^A\pi)\widetilde{\;}(a)(\xi)(s)=
X\dashind_B^A\pi(\alpha_s(a))(\xi(s))\ \mbox{ and }\
\rho_t(\xi)(s)=\Delta(t)^{1/2}\xi(st).
\]
The next lemma is similar to \cite[Theorem~1]{aHRqs}.

\begin{lemma}\label{subregularity}
Let $\pi$ be a faithful nondegenerate representation of $B$ on $\H_\pi$.
There is an isometry $W$ of $X_0\otimes_D\H_\pi$ into
$L^2(G,X\otimes_B\H_\pi)$ such that
\[
W(x\otimes_D \pi(b)h)(s)=\gamma_s(x)\otimes_B \pi(b)h \ \mbox{ for $x\in
X_0$, $b\in B_0$},
\]
and then $W(\nu,V)W^*$ is the restriction of the  regular
representation $\big((X\dashind_B^A\pi)\widetilde{\;},\rho\big)$ to the
range of $W$. In particular, $(\nu,V)$ is a covariant representation of
$(A,\alpha)$.
\end{lemma}

\begin{proof} 
We begin by noting that if $x\in X_0$ and $b\in B_0$, then
\begin{align*}
\int_G\|\gamma_s(x)\cdot b\|^2\,ds&=\int_G\|\langle\gamma_s(x)\cdot b\,,\,
\gamma_s(x)\cdot b\rangle_B\|\,ds\\
&=\int_G\|b^*\beta_s(\langle x\,,\,x\rangle_B)b\|\,ds\\
&\leq \|b\|\int_G\|b^*\beta_s(\langle x\,,\,x\rangle_B)\|\,ds\\
&=\|b\|\int_G \|\langle\gamma_s(x)\cdot b\,,\,\gamma_s(x)\rangle_B\|\,ds\\
&\leq\|b\|\,\|x\|\int_G \|\gamma_s(x)\cdot b\|\, ds,
\end{align*}
which is finite by Definition~\ref{defn-pr}(2).
Since $\|\gamma_s(x)\otimes_B\pi(b)h\|\leq\|\gamma_s(x)\cdot b\|\,\|h\|$
it follows that $W$
maps
$X_0\odot\pi(B_0)\H_\pi$ into $L^2(G,X\otimes_B\H_\pi)$. To see that
$W$ is isometric, we fix two vectors $x\otimes_D\pi(b)h$ and $y\otimes_D
\pi(c)k$ in $X_0\odot\pi(B_0)\H_\pi$, and compute:
\begin{align*}
\big(W(x\otimes_D \pi(b)h)\,\big|\, W(y\otimes_D \pi(c)k))\big)
&=\int_G \big(W(x\otimes_D \pi(b)h)(s)\,\big|\, W(y\otimes
\pi(c)k)(s)\big)\,ds\\ 
&=\int_G \big(\gamma_s(x)\otimes_B
\pi(b)h\,\big|\,\gamma_s(y)\otimes_B \pi(c)k\big)\,ds\\ 
&=\int_G \big(\pi
(c^*\langle\gamma_s(y)\,,\,\gamma_s(x)\rangle_Bb)h\,\big|\,k\big)\,ds\\
&=\bigg(\pi \Big(\int_G c^*\beta_s(\langle y\,,\,x\rangle_B)b\,ds\Big)
h\;\bigg|\;k\bigg)\\ 
&=\big(\pi(\langle y\,,\,x\rangle_D) \pi(b)h\,\big|\,
\pi(c)k\big)\\
&=\big(x\otimes_D \pi(b)h\,\big|\, y\otimes_D \pi(c)k\big).
\end{align*}
Thus $W$ extends to an isometry on
$X_0\otimes_D\H_\pi=\overline{X_0\odot\pi(B_0)\H_\pi}$.

Now for $x\otimes_Dh\in X_0\odot\pi(B_0)\H_\pi$, we have
\begin{align*}
\big(W\nu(a)(x\otimes_D h)\big)(s)
&= W(a\cdot x\otimes_D h)(s)=\gamma_s(a\cdot x)\otimes_B h\\
&=\alpha_s(a)\cdot\gamma_s(x)\otimes_B h
=X\dashind_B^A\pi(\alpha_s(a))\big(W(x\otimes_D h)(s)\big)\\
&=\big((X\dashind_B^A\pi)\widetilde{\;}(a)W(x\otimes_D h)\big)(s),
\end{align*}
which proves the first intertwining relation. To check that
$WV_s(x\otimes_D h)=\rho_sW(x\otimes_D h)$ is even easier.
\end{proof}

\begin{remark}
It is crucial in this argument that we start with a representation $\pi$
of $B$, and the lemma suggests that this choice is one of the reasons
that the reduced crossed product is appropriate on the left-hand side.
The importance of this issue is emphasized by Example~2.1 of
\cite{rie-pr}, where
$B=\K(L^2(G))$ and $D_0$ is contained in the subalgebra $\rho(C_c(G))$ of
$M(B)=B(L^2(G))$. Any representation $U$ of $G$ gives a representation of
$C_c(G)\cong\rho(C_c(G))$, which extends to a representation of
$\K(L^2(G))\cong C_0(G)\rtimes G$ precisely when there is a compatible
representation of $C_0(G)$; by the imprimitivity theorem, this happens
precisely when
$U$ is induced from a representation of $\{e\}$, and hence is a regular
representation of
$G$.
\end{remark}

\begin{lemma}\label{actionofE}
For $x,y\in X_0$, define $_E\langle x,y\rangle:G\to A_0$ by
\[
_E\langle x,y\rangle(s)=\Delta(s)^{-1/2}{}_A\langle x\,,\,\gamma_s(y)\rangle,
\]
which belongs to $L^1(G,A)$ by Definition~\ref{defn-pr}(1). Let $\theta_{x,y}$ be the compact operator on
$(X_0)_D$ defined by $\theta_{x,y}z=x\cdot\langle y\,,\,z\rangle_D$. Then
\begin{equation}\label{EactsonX}
\nu\rtimes V({}_E\langle x\,,\,y\rangle)=X_0\dashind_D^\L\bar\pi(\theta_{x,y}).
\end{equation}
\end{lemma}

\begin{proof}
For $z\in X_0$, $h\in \H_\pi$ and $c\in B_0$, we have
\begin{align}
\nu\rtimes V({}_E\langle x,y\rangle)(z\otimes_D\pi(c)h)
&=\int_G \nu(\Delta(s)^{-1/2}{}_A\langle
x\,,\,\gamma_s(y)\rangle)\Delta(s)^{1/2}(\gamma_s(z)\otimes_D\pi(c)h)\,ds\notag\\
&=\int_G \big({}_A\langle
x\,,\,\gamma_s(y)\rangle\cdot\gamma_s(z)\big)\otimes_D\pi(c)h\,ds\notag\\
&=\int_G \big(x\cdot\langle
\gamma_s(y)\,,\,\gamma_s(z)\rangle_B\big)\otimes_D\pi(c)h\,ds\notag\\
&=\int_G \big(x\cdot\beta_s(\langle
y\,,\,z\rangle_B)\big)\otimes_D\pi(c)h\,ds.\label{last}
\end{align}
At this stage we want to whip the $ds$ past $\otimes_D\pi(c)h$: then
the integral would be that defining $x\cdot\langle
y\,,\,z\rangle_D=\theta_{x,y}(z)$, and we'd be done. Unfortunately, the
resulting integral converges in the norm coming from the $B$-valued inner
product, so we have to work to pull it through a balanced tensor product
defined using the $D$-valued inner product. 

The
$(X_0\otimes_D\H_\pi)$-valued integral in \eqref{last} is characterized by
its inner products with vectors of the form $w\otimes_D\pi(b)k$ for $b\in
B_0$:
\begin{align}
\Big(\int_G \big(x\cdot\beta_s(\langle
y\,,\,z\rangle_B)\big)&\otimes_D\pi(c)h\,ds\;\Big|
\;w\otimes_D\pi(b)k\Big)\notag\\
&=\int_G \Big(\big(x\cdot\beta_s(\langle
y\,,\,z\rangle_B)\big)\otimes_D\pi(c)h\;\Big|\;w\otimes_D\pi(b)k\Big)\,ds\notag\\
&=\int_G \Big(\pi\big(b^*\big\langle w\,,\,x\cdot\beta_s(\langle
y\,,\,z\rangle_B)\big\rangle_Dc\big)h\;\Big|\;k\Big)\,ds\notag\\
&=\int_G \bigg(\pi\Big(\int_G b^*\beta_t\big(\big\langle
w\,,\,x\cdot\beta_s(\langle
y\,,\,z\rangle_B)\big\rangle_B\big)c\,dt\Big)h\;\bigg|\;k\bigg)\,ds,\notag\\
\intertext{which, because the inner integral is that of a $B$-valued
function, is just}
&=\int_G\int_G  \Big(\pi\big(b^*\beta_t\big(\big\langle
w\,,\,x\cdot\beta_s(\langle
y\,,\,z\rangle_B)\big\rangle_B\big)c\big)h\;\Big|\;k\Big)\,dt\,ds\notag\\
&=\int_G\int_G  \Big(\pi\big(b^*\beta_t(\langle
w\,,\,x\rangle_B)\beta_{ts}(\langle
y\,,\,z\rangle_B)c\big)h\;\Big|\;k\Big)\,dt\,ds.\label{lastagain}
\end{align}
The two elements $b$ and $c$ of $B_0$ are there to ensure that the
integrand in this double integral is integrable on $G\times G$, so that
we can apply Fubini's Theorem to continue:
\begin{align}
\eqref{lastagain}
&=\int_G\int_G  \Big(\pi\big(b^*\beta_t(\langle
w\,,\,x\rangle_B)\beta_{ts}(\langle
y\,,\,z\rangle_B)c\big)h\;\Big|\;k\Big)\,ds\,dt\notag\\
&=\int_G\int_G  \Big(\pi\big(b^*\beta_t(\langle
w\,,\,x\rangle_B\beta_{s}(\langle
y\,,\,z\rangle_B)c)\big)h\;\Big|\;k\Big)\,ds\,dt\notag\\
&=\int_G\int_G  \Big(\pi\big(b^*\beta_t(\langle
w\,,\,x\rangle_B \beta_{s}(\langle
y\,,\,z\rangle_B)c)\big)h\;\Big|\;k\Big)\,dt\,ds.\label{andagain}
\end{align}
We can now go backwards through the previous analysis to see that
\[
\eqref{andagain}=\big(\theta_{x,y}z\otimes_D\pi(c)h\,\big|\,
w\otimes_D\pi(b)k\big)
=\big(X_0\dashind_D^\L\bar\pi(\theta_{x,y})(z\otimes_D\pi(c)h)\,\big|\,
w\otimes_D\pi(b)k\big),
\]
and the result follows.
\end{proof}

To get our Morita equivalence, we consider 
\[
E_0:=\sp\{{}_E\langle x\,,\,y\rangle:x,y\in X_0\}\subset L^1(G,A).
\]
Because we know from Lemma~\ref{subregularity} that the representation $(\nu,V)$ is equivalent to a
subrepresentation of the regular representation,
$\nu\rtimes V$ extends to a representation of the reduced crossed product
$A\rtimes_{\alpha,r}G$. Equation~\eqref{EactsonX} therefore implies that
$\nu\rtimes V$ carries the closure $E$ of $E_0$ in $A\rtimes_{\alpha,r}G$
into the image
$X_0\dashind_D^\L\bar\pi\big(\K\big(\overline{(X_0)_D}\big)\big)$ of
the imprimitivity algebra
\[
\K\big(\overline{(X_0)_D}\big):=\clsp\{\theta_{x,y}:x,y\in X_0\}\subset
\L\big(\overline{(X_0)_D}\big).
\]  
We would \emph{like} to prove next that $\nu\rtimes V$ is isometric for
the reduced norm on $E_0$. If we could do this, we could deduce that
$\nu\rtimes V$ is an isometric linear isomorphism of
$E$ onto $\K\big(\overline{(X_0)_D}\big)$; this would imply both that $E$
is a $C^*$-algebra (because $\K$ is, and $\nu\rtimes V$ is a $*$-algebra
homomorphism on $A\rtimes_{\alpha,r}G$), and that $\overline{(X_0)_D}$ is
an $E$--$D$ imprimitivity bimodule.

The program of the previous paragraph works without any problems when the
isometry $W$ of Lemma~\ref{subregularity} maps
$X_0\otimes_D\H_\pi$ onto
$L^2(G,X\otimes_B\H_\pi)$, or, equivalently, when 
\begin{equation}\label{rangeW}
\sp\{s\mapsto \gamma_s(x)\otimes_D\pi(b)h: x\in X_0,\;b\in B_0,\;h\in
\H_\pi\}
\end{equation}
is dense in $L^2(G,X\otimes_B\H_\pi)$. To make this independent of the
choice of Hilbert space, we could ask instead that the corresponding map
of
$X_0\otimes_D B$ into
$L^2(G,X_B)$ be an isomorphism of Hilbert $B$-modules. In the
examples  where $X_0$, $A_0$ and $B_0$ consist of
functions of compact support, this seems remarkably similar to asking that
the maps
$s\mapsto {}_E\langle x\,,\,y\rangle(s)$ span a dense subspace of $L^1(G,A)$.
So what we have proved at this stage is already potentially interesting:

\begin{prop}
Suppose that the functions $\{s\mapsto \gamma_s(x)\cdot b:x\in X_0,b\in
B_0\}$ span a dense subspace of $L^2(G,X_B)$, and that
the space
$E_0$ is dense in
$L^1(G,A)$. Then $\overline{(X_0)_D}$ is an $A\rtimes_{\alpha,r} G$--$D$
imprimitivity bimodule.
\end{prop}

Unfortunately, for arbitrary proper actions we do not see how to prove
directly that $\nu\rtimes V$ is isometric on $E$. So to establish the
Morita equivalence of $E$ and $D$, we prove that $X_0$ is an
$E_0$--$D_0$ pre-imprimitivity bimodule. Since we already know that $X_0$
is a pre-Hilbert $D$-module, it remains to show that
$E_0$ is a $*$-algebra which acts by bounded operators on $(X_0)_D$
according to the formula ${}_E\langle x\,,\,y\rangle\cdot z=x\cdot\langle
y\,,\,z\rangle_D$, that ${}_E\langle\cdot\,,\,\cdot\rangle$ is then a pre-inner
product with respect to the completion $E$ of $E_0$ (that is, the closure
$E$ of
$E_0$ in the reduced norm), and that $D$ acts by bounded operators on
${}_E(X_0)$.

\begin{lemma}\label{formulas}
The pairing ${}_E\langle\cdot\,,\,\cdot\rangle$ on $X_0$ satisfies
\[
{}_E\langle x\,,\,y\rangle* {}_E\langle z,w\rangle
={}_E\langle x\cdot\langle y\,,\, z\rangle_D,w\rangle
\ \mbox{ and }\ {}_E\langle x\,,\,y\rangle^*={}_E\langle y\,,\,x\rangle,
\]
where the product and adjoint are those of $L^1(G,A)\subset
A\rtimes_\alpha G$.
\end{lemma} 

\begin{proof}
For the first identity, we compute
\begin{align*}
{}_E\langle x\,,\, y\rangle*{}_E\langle z\,,\, w\rangle(s)
&=\int_G \Delta(r)^{-1/2}{}_A\langle x\,,\,\gamma_r(y)\rangle
\alpha_r\big(\Delta(r^{-1}s)^{-1/2}{}_A\langle z\,,\,
\gamma_{r^{-1}s}(w)\rangle\big)\, dr\\
&=\int_G {}_A\langle x\,,\,\gamma_r(y)\rangle {}_A\langle
\gamma_r(z)\,,\,\gamma_s(w)\rangle\, dr\,\Delta(s)^{-1/2}\\
&=\int_G {\phantom{\big|}}_A\big\langle {}_A\langle x\,,\,\gamma_r(y)\rangle
\cdot\gamma_r(z)\,,\,\gamma_s(w)\big\rangle\, dr\, \Delta(s)^{-1/2}\\
&={\phantom{\Big|}}_A\Big\langle\int_G x\cdot \langle \gamma_r(y)\,,\,
\gamma_r(z)\rangle_B\, dr\,,\;\gamma_s(w)\Big\rangle\;\Delta(s)^{-1/2}\\
&={\phantom{\big|}}_A\big\langle x\cdot\langle y\,,\,
z\rangle_D\,,\,\gamma_s(w)\big\rangle\,\Delta(s)^{-1/2}\\
&={\phantom{\big|}}_E\big\langle x\cdot \langle y\,,\, z\rangle_D\,,\,
w\big\rangle(s).
\end{align*}
The second identity follows from a simple algebraic manipulation.
\end{proof}

\begin{prop}
The set $E_0:=\sp\{{}_E\langle x\,,\,y\rangle:x,y\in X_0\}$ is a
$*$-subalgebra of $L^1(G,A)$, and there is a left action of $E_0$ on
$(X_0)_D$ such that  
\[
{}_E\langle x\,,\,y\rangle\cdot z=x\cdot\langle
y\,,\,z\rangle_D \ \mbox{ and }\ \langle e\cdot x,e\cdot x\rangle_D \leq
\|e\|^2\langle x\,,\,x\rangle_D\mbox{ as elements of $D$,}
\]
for $x,y,z\in X_0$ and $e\in E_0$.
\end{prop}

\begin{proof}
Because $x\cdot\langle
y\,,\,z\rangle_D$ belongs to $X_0$, the formulas in Lemma~\ref{formulas} show
that $E_0$ is a $*$-subalgebra of $L^1(G,A)$. We know from Lemma~\ref{actionofE} that
the $*$-homomorphism $\phi:=\big(\Ind_D^\L\bar\pi\big)^{-1}\circ
(\nu\rtimes V)$ restricts to a $*$-homomorphism of $E_0$ into
$\L\big(\overline{(X_0)_D}\big)$ such that $\phi({}_E\langle
x\,,\,y\rangle)=\theta_{x,y}$, and $\phi$ gives the required action of $E_0$
on
$(X_0)_D$: $e\cdot z$ is by definition $\phi(e)z$. For any bounded
operator
$T\in
\L\big(\overline{(X_0)_D}\big)$, we have $\langle Tx\,,\,Tx\rangle_D \leq
\|T\|^2\langle x\,,\,x\rangle_D$; because $\Ind_D^\L\bar\pi$ is isometric and
$\nu\rtimes V$ is decreasing for the reduced norm, we have
$\|\phi(e)\|\le \|e\|$, and the inequality follows. 
\end{proof}

\begin{prop}
The pairing ${}_E\langle\cdot\,,\,\cdot\rangle$ is an inner product with
values in the $C^*$-algebra $E:=\overline{E_0}$, and 
\begin{equation}\label{DactsonEX}
{}_E\langle x\cdot d\,,\,x\cdot d\rangle\leq \|d\|^2
{}_E\langle x\,,\,x\rangle\ \mbox{ as elements of the $C^*$-algebra $E$}.
\end{equation}
\end{prop}

\begin{proof}
We know from Lemma~\ref{formulas} that
\[
{}_E\langle x\,,\,y\rangle* {}_E\langle z,w\rangle
={}_E\langle x\cdot\langle y\,,\, z\rangle_D,w\rangle;
\]
since the left action of $E_0$ satisfies ${}_E\langle
x\,,\,y\rangle\cdot z=x\cdot\langle y\,,\, z\rangle_D$, it follows that
$e\cdot{}_E\langle z\,,\,w\rangle={}_E\langle e\cdot z\,,\,w\rangle$ for all
$e\in E_0$. Lemma~\ref{formulas} also shows that ${}_E\langle
x\,,\,y\rangle^*={}_E\langle y\,,\,x\rangle$, so ${}_E\langle\cdot\,,\,\cdot\rangle$
has the required algebraic properties. 

To see positivity, we fix a representation $\rho$ of $A$ on $\H$, and
consider the left-regular representation $(\widetilde\rho,\lambda)$ on
$L^2(G,\H)$ given by
\[
(\widetilde\rho(a)\xi)(s)=\rho(\alpha_{s^{-1}}(a))\xi(s)\ \mbox{ and }\ 
(\lambda_t\xi)(s)=\lambda(t^{-1}s).
\]
Now we let $x\in X_0$ and consider $\xi\in L^2(G,\H)$ of the form
$s\mapsto f(s)h$ for $f\in C_c(G)$. Then we have, by a variant of Remark~\ref{pullpinbimods},
\begin{align*}
\big(\widetilde\rho\rtimes\lambda( {}_E\langle x\,,\, x\rangle
)\xi\,\big|\,\xi \big) 
&= \Big(\int_G\tilde\rho\big( {}_E\langle x\,,\,x\rangle
(s)\big)
\lambda_s \xi\, ds\;\Big|\; \xi\Big)\\
&=\int_G\bigg(\Big(\int_G\tilde\rho\big(\Delta(s)^{-1/2}{}_A\langle
x\,,\,
\gamma_s(x)\rangle\lambda_s\xi\, ds\Big)(t)\;\bigg|\;\xi(t)\bigg)\,dt\\
&=\int_G\!\int_G \Big(\rho\big(\alpha_{t^{-1}}({}_A\langle
x\,,\,\gamma_s(x)\rangle)\big)(\xi(s^{-1}t))\;
\Big|\;\xi(t)\Big)\Delta(s)^{-1/2}\,
ds\, dt\\
&=\int_G\!\int_G \Big(\rho\big({}_A\langle
\gamma_{t^{-1}}(x)f(t)\,,\,\gamma_{t^{-1}s}(x)f(s^{-1}t)\rangle\big)h\;
\Big|\;h\Big)\Delta(s)^{-1/2}\,
ds\, dt.
\end{align*}
Since the function $t\mapsto \gamma_{t^{-1}}(x)f(t)$ belongs to
$C_c(G,X)$, we can apply Fubini's Theorem, and then substitute
$r=s^{-1}t$ to reduce this to
\[
\int_G\!\int_G \Big(\rho\big({\phantom{\big|}}_A\big\langle
\gamma_{t^{-1}}(x)f(t)\Delta(t)^{-1/2}\,,\,
\gamma_{r^{-1}}(x)f(r)\Delta(r)^{-1/2}\big\rangle\big)h\;
\Big|\;h\Big)\,
dr\, dt,
\]
which has the form $\big(\rho({}_A\langle y\,,\,y\rangle)h\,\big|\,h\big)$
for the element $y$ of $X$
defined by 
\[
y:=\int_G\gamma_{t^{-1}}(x)f(t)\Delta(t)^{-1/2}\,dt,
\] 
and hence
is positive. Since $\xi$ of the given form span a dense subspace of
$L^2(G,\H)$, this proves that 
$\widetilde\rho\rtimes\lambda( {}_E\langle x\,,\, x\rangle
)$ is positive. 

Thus ${}_E\langle \cdot\,,\, \cdot\rangle$ is a pre-inner product; it is
definite because the regular representation is faithful on the reduced crossed product where $E$ sits. If now $d\in
D\subset M(B)^\beta$, then $\gamma_s(x\cdot
d)=\gamma_s(x)\cdot\beta_s(d)=\gamma_s(x)\cdot d$, so we can repeat the
calculation of the previous paragraph to see that
\[
\big(\widetilde\rho\rtimes\lambda(\|d\|^2
{}_E\langle x\,,\,x\rangle-{}_E\langle x\cdot d\,,\,x\cdot
d\rangle)\xi\,\big|\,\xi\big)=
\big(\rho(\|d\|^2
{}_A\langle y\,,\,y\rangle-{}_A\langle y\cdot d\,,\,y\cdot
d\rangle)h\,\big|\,h\big),
\]
which is positive because $M(B)$ acts as bounded operators on
${}_AX$. This proves (\ref{DactsonEX}).
\end{proof}

\begin{remark}
As in \cite{rie-pr}, $E$ is an ideal in $A\rtimes_{\alpha,r}G$. To see
this, one observes that the dense subalgebra $A_0$ of $A$ satisfies
$A_0E_0\subset E_0$, that $E_0$ is $\alpha$-invariant, and deduce that
$fE_0\subset E=\overline{E_0}$ for any $f$ in the dense subalgebra
$C_c(G,A_0)$ of $A\rtimes_{\alpha,r}G$.
\end{remark}

\begin{definition} Following \cite{rie-pr}, we say that a proper action $\gamma$ on ${}_AX_B$ is \emph{saturated} with respect to $X_0$ if $E=A\rtimes_{\alpha,r}G$.
\end{definition}

To sum up, we have now proved:

\begin{thm}\label{thm-sumup}
Suppose the Morita equivalence ${}_{(A,\alpha)}(X,\gamma)_{(B,\beta)}$ is
proper in the sense of Definition~\ref{defn-pr} with respect to the
pre-imprimitivity bimodule ${}_{A_0}(X_0)_{B_0}$. Then the completion of
$X_0$ in the norm $\|x\|:=\|\langle x\,,\,x\rangle_D\|^{1/2}$ implements a
Morita equivalence between the ideal
\[
E:=\clsp\{s\mapsto \Delta(s)^{-1/2}{}_A\langle
x\,,\,\gamma_s(y)\rangle:x,y\in X_0\}\subset A\rtimes_{\alpha,r}G
\]
and the $C^*$-algebra
\[
D:=\clsp\{\langle x\,,\,y\rangle_D:x,y\in X_0\}\subset M(B).
\]
In particular, if the action $\gamma$ is saturated, then $\overline{(X_0)_D}$
implements a Morita equivalence between $A\rtimes_{\alpha,r}G$ and $D$.
\end{thm}

The general theory produces a pre-imprimitivity bimodule on which only the spans of the ranges of the inner products act.  In the main  examples, there are algebras of continuous functions of compact support which ought to act too, and it is important that the formulas extend. The following lemma gives conditions under which   the extended left action is given by the expected formula.

\begin{lemma}\label{formulaforaction}
Let $f:G\to A$ be a  continuous function such that both
$s\mapsto \|f(s)\|$ and $s\mapsto \|f(s)\|\Delta(s)^{1/2}$ are integrable.
Suppose that $x\in X_0$ and that the integral $\int_G
f(s)\cdot\gamma_s(x)\Delta(s)^{1/2}\,ds$, which converges in $X_B$
because of our second integrability hypothesis on $f$, belongs to $X_0$.
Let
$\pi$ be a representation of
$B$ and let $(\nu,V)$ be the covariant representation discussed in
Lemma~\ref{subregularity}; note that the first integrability
hypothesis on $f$ implies that $\nu\rtimes V(f)$ makes sense as a bounded
operator on
$X_0\otimes_D\H_\pi$. Then
\[
\big({\textstyle \int_G
f(s)\cdot\gamma_s(x)\Delta(s)^{1/2}\,ds}\big)\otimes_D h=\nu\rtimes
V(f)(x\otimes_D h) \ \mbox{ for }h\in \H_\pi.
\]
In other words, the left action of $f\in
L^1(G,A)$ on $x\in (X_0)_D$ is given by
the integral formula
\[
f\cdot x=\int_G
f(s)\cdot\gamma_s(x)\Delta(s)^{1/2}\,ds.
\]
\end{lemma}

\begin{proof}
To make things a little easier on the eye, we shall write 
$g(s):=f(s)\cdot\gamma_s(x)\Delta(s)^{1/2}$.
We fix $y\otimes_D\pi(b)k\in X_0\otimes \pi(B_0)\H_\pi$, and compute:
\begin{align}
\big(\nu\rtimes
V(f)(x\otimes_D h)\,\big|\,y\otimes\pi(b)k\big)
&=\Big(\int_G\nu(f(s))V_s(x\otimes_D
h)\,ds\,\big|\,y\otimes\pi(b)k\Big)\notag\\
&=\Big(\int_G(g(s)\otimes_D
h)\,ds\,\big|\,y\otimes\pi(b)k\Big)\notag\\
&=\int_G (g(s)\otimes_D h\,|\,y\otimes\pi(b)k)\,ds\notag\\
&=\int_G \big(\pi(b^*\langle
y\,,\,g(s)\rangle_D)h\,\big|\,k\big)\,ds\notag\\
&=\int_G \Big(\pi\Big(\int_G b^*\beta_t(\langle
y\,,\,g(s)\rangle_B)\,dt\Big)h\,\big|\,k\Big)\,ds\notag\\
&=\int_G\!\int_G \big(\pi(b^*\beta_t(\langle
y\,,\,g(s)\rangle_B))h\,\big|\,k\big)\,dt\,ds,
\label{maybe?}
\end{align}
using standard properties of $B$-valued integrals.
On the other hand, 
\begin{align*}
\Big(\Big(\int_G g(s)\,ds\Big)\otimes_D
h\,\big|\,y\otimes\pi(b)k\Big) 
&=\Big(\pi\Big(b^*\big\langle
y\,,\, \int_G g(s)\,ds\big\rangle_D\Big) h\,\big|\,k\Big)\notag\\
&=\Big(\pi\Big(\int_G b^*\beta_t(\langle
y\,,\, \int_G g(s)\,ds\rangle_B)\,dt\Big) h\,\big|\,k\Big)\notag\\
&=\int_G \Big(\pi\Big( b^*\beta_t(\langle
y\,,\, \int_G g(s)\,ds\rangle_B)\Big)
h\,\big|\,k\Big)\,dt.
\end{align*}
Because the inside integral converges in $X_B$, we can pull it
through the $B$-valued inner product with $y$; now we have an ordinary
$B$-valued integral, and we can pull the automorphisms and
representation through to recover
\begin{equation}\label{atlast?}
\int_G\!\int_G \big(\pi(b^*\beta_t(\langle
y\,,\,g(s)\rangle_B))h\,\big|\,k\big)\,ds\,dt.
\end{equation}
But now we're talking about ordinary scalar-valued integrals; the
element $b$ is a sum of elements of the form $\langle w\,,\,z\rangle_B$,
and for such an element the integrand
\[
\big(\pi(b^*\beta_t(\langle
y\,,\,g(s)\rangle_B))h\,\big|\,k\big)=
\big(\pi(\langle{}_A\langle \gamma_t(y)\,,\,w\rangle\cdot z\,,\,\gamma_t(g(s))
\rangle_B))h\,\big|\,k\big)
\] 
is integrable on $G\times G$. Thus an application of Fubini's Theorem
identifies \eqref{atlast?} with \eqref{maybe?}, and we are done.
\end{proof}

\section{Tensor-product decompositions of imprimitivity bimodules.}\label{decomp}

Suppose the Morita equivalence ${}_{(A,\alpha)}(X,\gamma)_{(B,\beta)}$ is
proper  with respect to the
pre-imprimitivity module ${}_{A_0}(X_0)_{B_0}$, so that Theorem~\ref{thm-sumup} gives a Morita equivalence $\overline{(X_0)}$ between an ideal in $A\rtimes_{\alpha,r} G$ and a generalized fixed-point algebra $D$ for $(B, \beta)$. We now want to show that $(B,\beta)$ is itself proper in Rieffel's sense \cite[Definition~1.2]{rie-pr}, and to relate $\overline{(X_0)}$ to other Morita equivalences involving Rieffel's generalized fixed-point algebra $B^\beta$. Both $D$ and $B^\beta$ are by definition subalgebras of $M(B)^\beta$, but it is not obvious that they must be the same subalgebra. Indeed, it is not even obvious that we get the same generalized fixed-point algebra when $\beta$ is proper with respect to two different dense $*$-subalgebras. Fortunately, we have been able to show that at least when $\gamma$ is saturated, all the  fixed-point algebras relevant to us coincide. When we have sorted this out, it will be relatively easy to get the desired relations between imprimitivity bimodules. 

We begin by giving a careful statement of our main results. Since we are concerned about possibly different fixed-point algebras, we shall denote by $(B,B_1)^\beta$ the generalized fixed-point algebra as defined in \cite{rie-pr} when $\beta$ is proper with respect to a particular subalgebra $B_1$.

\begin{thm}\label{thm-decomposition}
Suppose that the Morita equivalence ${}_{(A,\alpha)}(X,\gamma)_{(B,\beta)}$ is proper  with respect to ${}_{A_0}(X_0)_{B_0}$, with generalized fixed-point algebra $D$. Then
\begin{enumerate}
\item[(1)] the action $\beta$ on $B$ is proper with respect to the subalgebra
\[
B_1:=\langle X_0\,,\,X_0\rangle_B:=\sp\{\langle x\,,\,y\rangle_B:x,y\in X_0\}, 
\]
and the generalized fixed-point algebra $(B,B_1)^\beta$ is an ideal in $D$; 
\item[(2)] the action $\gamma$ is also proper with respect to the smaller pre-imprimitivity bimodule ${}_{A_0}(X_0)_{B_1}$, and then has the same generalized fixed-point algebra $D$.
\end{enumerate}
The action  $\beta$ is saturated with respect to $B_1$ if and only if $\gamma$ is saturated with respect to $X_0$, and then $(B,B_1)^\beta=D$. 
There is an isomorphism 
\begin{equation*}
\Omega:\big( X\rtimes_r G\big)\otimes_{B\rtimes_r G} \overline{B_1}\to\overline{X_0}
\end{equation*}
of $(A\rtimes_{\alpha,r} G)$--$D$-imprimitivity bimodules such that
\begin{equation}\label{eq-isom}
\Omega(f\otimes b)= \int_G f(s)\cdot\beta_s(b)\Delta(s)^{1/2}\, ds
\end{equation}
where $b\in B_1$ and $f$ has the form $f(s):=x\cdot\beta_s(c^*)\Delta(s)^{-1/2}$ with $x\in X_0$ and $c\in B_1$.
\end{thm}

\begin{proof}[Proofs of \textnormal{(1)} and \textnormal{(2)}]
For (1), we need to verify the two items of \cite[Definition~1.2]{rie-pr}. Write $F$ for $(B,B_1)^\beta$, and let $b=\langle v\,,\,w\rangle_B$ and $c=\langle x\,,\,y\rangle_B$ be typical spanning elements of $B_1$. Then
\[
b\beta_s(c^*)
=\big\langle v\,,\, w\cdot\beta_s(\langle y\,,\, x\rangle_B)\big\rangle_B 
=\big\langle v\,,\,{}_A\langle w\,,\,\gamma_s(y)\rangle\cdot\gamma_s(x)\big\rangle_B,
\]
and hence 
\[
\|b\beta_s(c^*)\|\leq \|v\|\|x\|\|{}_A\langle w\,,\,\gamma_s(y)\rangle\|.
\]
The function $s\mapsto {}_A\langle w\,,\, \gamma_s(y)\rangle$ and its product with $s\mapsto \Delta(s)^{1/2}$ are in $L^1(G,A)$ because $\gamma$ is proper, so it follows that $s\mapsto b\beta_s(c)$ and $s\mapsto \Delta(s)^{1/2}b\beta_s(c)$ are integrable. This gives the first item of \cite[Definition~1.2]{rie-pr}.

Set $z:=v\cdot\langle x\,,\,y\rangle_{B}$. Then $z\in X_0$, and  
\begin{equation}
\label{eq42}
b^*c
=\langle w\,,\, v\rangle_B\langle x\,,\, y\rangle_B
=\big\langle w\,,\, v\cdot\langle x\,,\, y\rangle_B\big\rangle_{B}
=\langle w\,,\,z\rangle_B;
\end{equation}
Definition~\ref{defn-pr}(3) says there is a multiplier $\langle w\,,\, z\rangle_D\in M(B_0)^\beta$ such that 
\begin{equation}\label{firstbit}
\int_G a\beta_s(b^*c)\, ds\stackrel{\eqref{eq42}}{=}\int_G a\beta_s(\langle w\,,\, z\rangle_B)\, ds =a\cdot\langle w\,,\, z\rangle_D\ \mbox{for all $a\in B_1\subset B_0$}.
\end{equation}
We claim that $\langle w\,,\, z\rangle_{D}$ multiplies $B_1$ (we already know it multiplies $B_0$).
If  $b'=\langle x_1\,,\, x_2\rangle_B\in B_1$, then $b'\cdot\langle w\,,\,z\rangle_D=\big\langle x_1\,,\, x_2\cdot\langle w\,,\, z\rangle_D\big\rangle_B$ because $\langle w\,,\, z\rangle_D\in M(B)$. But $x_2\cdot\langle w\,,\, z\rangle_D$ is back in $X_0$, so $b'\cdot\langle w\,,\,z\rangle_D\in\langle X_0\,,\, X_0\rangle_B=B_1$. Thus $\langle w\,,\, z\rangle_{D}\in M(B_1)$. 

We now define $\langle b\,,\, c\rangle_{F}:=\langle w\,,\, z\rangle_D$, and \eqref{firstbit} gives the second item of \cite[Definition~1.2]{rie-pr}. Notice that by definition of $\langle\cdot\,,\,\cdot\rangle_F$, we have $F=\overline{\langle X_0\,,\, X_0\cdot B_1\rangle_D}$.
 Since  $X_0\cdot B_1$ is a sub-module of $_{E_0}(X_0)_{D_0}$, it follows from the Rieffel correspondence  that $F$ is an ideal of $D$. This gives (1).

For (2), we have to verify the three properties of Definition~\ref{defn-pr} for ${}_{A_0}(X_0)_{B_1}$. Since $B_1\subset B_0$, the integrability properties are clear. So it suffices to check Definition~\ref{defn-pr}(3) and to show that $D$ and $D_1:=D({}_{A_0}(X_0)_{B_1})$ coincide.
If $x,y,z\in X_0$ then $z\cdot \langle x\,,\, y\rangle_{D}$ is in $X_0$, so 
\[
\langle w\,,\, z\rangle_B\langle x\,,\, y\rangle_{D}=\langle w\,,\, z\cdot \langle x\,,\, y\rangle_{D}\rangle_B
\]
is in $B_1$. Thus $\langle x\,,\, y\rangle_{D}$ multiplies $B_1$, and $\langle x\,,\, y\rangle_{D_1}:=\langle x\,,\, y\rangle_{D}$ has the properties described in Definition~\ref{defn-pr}(3). But with this definition of $\langle \cdot\,,\, \cdot\rangle_{D_1}$ we trivially have $D=D_1$.
\end{proof}

The proof of the decomposition isomorphism \eqref{eq-isom} uses a general lemma about imprimitivity bimodules over a linking algebra. If ${}_AX_B$ is an imprimitivity bimodule, let $\widetilde X$ be the conjugate vector space, so that there is an additive bijection $\flat:X\to\widetilde X$ such that $\flat(\lambda\cdot x)=\bar\lambda\flat(x)$. Then $\widetilde X$ is a $B$--$A$ imprimitivity bimodule with
\begin{gather*}
b\cdot \flat(x)=\flat(x\cdot b^*)\quad
\flat(x)\cdot a=\flat(a^*\cdot x)\\
{}_B\langle\flat(x)\,,\,\flat(y)\rangle=\langle x\,,\, y\rangle_B\quad
\langle\flat(x)\,,\,\flat(y)\rangle_A={}_A\langle x\,,\, y\rangle.
\end{gather*}
Recall that the \emph{linking algebra} $L(X)$ of an imprimitivity bimodule ${}_AX_B$ is the collection of $2\times 2$ matrices
\[
L(X):=\Big\lbrace \begin{pmatrix}a&x\\\flat(y)&b\end{pmatrix}:a\in A, b\in B, x,y\in X\Big\rbrace
\]
with multiplication and involution given by
\[
\begin{pmatrix}a&x\\\flat(y)&b
\end{pmatrix}
\begin{pmatrix}a'&x'\\\flat(y')&b'
\end{pmatrix}
:=\begin{pmatrix}
aa'+{}_A\langle x\,,\,y'\rangle 
&a\cdot x'+x\cdot b'\\
\flat(y)\cdot a'+b\cdot\flat(y')
&\langle y\,,\, x'\rangle_B + bb'
\end{pmatrix}
\]
and
\[\begin{pmatrix}a&x\\\flat(y)&b
\end{pmatrix}^*:=\begin{pmatrix}a^*&y\\\flat(x)&b^*
\end{pmatrix};
\]
$L(X)$ acquires a $C^*$-algebra structure by identifying it with the $C^*$-algebra $\K(X\oplus B)$ of compact operators on the right-Hilbert module direct sum $(X\oplus B)_B$ (see \cite[Corollary~3.21]{tfb}). The matrices 
\[
p=p_{L(X)}=\begin{pmatrix}1_{M(A)}&0\\0&0\end{pmatrix}\text{\ and } q=q_{L(X)}=\begin{pmatrix}0&0\\0&1_{M(B)}\end{pmatrix}
\]
define full projections in $M(L(X))$, and the maps
\begin{equation}\label{idsinlinking}
a\mapsto \begin{pmatrix}a&0\\0&0\end{pmatrix},\quad
b\mapsto \begin{pmatrix}0&0\\0&b\end{pmatrix}\quad\text{and}\quad
x\mapsto \begin{pmatrix}0&x\\0&0\end{pmatrix}
\end{equation}
identify $A$, $B$ and $X$ with the corners $pL(X)p$, $qL(X)q$ and $pL(X)q$ in $L(X)$, respectively.

The next lemma is a variation of \cite[Lemma~4.6]{es} and \cite[Proposition~4.3]{EKQR}, where the algebra $C$ is also required to be a linking algebra. In it, we use the identifications \eqref{idsinlinking} to produce actions of $A$, $B$ and $X$ on a module over $L(X)$. 

\begin{lemma}\label{lem-ibm-isom}
Let $X$ be an $A$--$B$ imprimitivity bimodule with linking algebra $L(X)$.
If $Z$ is an $L(X)$--$C$ imprimitivity
bimodule, then 
$pZ$  and $qZ$ are $A$--$C$ and $B$--$C$ imprimitivity bimodules, respectively, and there is an   isomorphism $\Omega:X\otimes_B qZ\to pZ$ of $A$--$C$ imprimitivity bimodules such that
\begin{equation}\label{keyiso}
\Omega(x\otimes_B qz)=x\cdot qz.
\end{equation}
\end{lemma}

\begin{proof}
Since $A=pL(X)p\subset L(X)$, it is easy to see that $pZ$ is an $A$-module; on $pZ$, the $L(X)$-valued inner product takes values in $pL(X)p$, and with ${}_A\langle pz\,,\, pz'\rangle:=p{}_{L(X)}\langle z\,,\,z'\rangle p$, $pZ$ becomes a full left Hilbert $A$-module. The right actions and inner products are already defined; the only thing we need to worry about is whether $pZ$ is full as a Hilbert $C$-module. So let $I$ be the ideal in $C$ spanned by the elements $\langle pz\,,\, pz'\rangle_C$. Then
\begin{align*}
Z\dashind_C^{L(X)}(I)&=\clsp\lbrace {}_{L(X)}\langle  z\cdot i\,,\, z'\rangle : z,z'\in Z, i\in
I\rbrace\\
&=\clsp\lbrace {}_{L(X)}\langle z\cdot \langle pw\,,\,
pw'\rangle_C\,,\, z'\rangle : z,z',w,w'\in Z\rbrace\\
&=\clsp \lbrace{}_{L(X)}\langle {}_{L(X)}\langle z\,,\, pw\rangle\cdot pw'\,,\,
z'\rangle : z,z',w,w'\in Z\rbrace\\
&=\clsp \lbrace {}_{L(X)}\langle z\,,\, pw\rangle {}_{L(X)}\langle
pw'\,,\, z'\rangle :z,z',w,w'\in Z\rbrace\\
&=\clsp \lbrace {}_{L(X)}\langle z\,,\, w\rangle p^*p{}_{L(X)}\langle
w'\,,\, z'\rangle :z,z',w,w'\in Z\rbrace\\
&=\overline{L(X)pL(X)},\end{align*}
which is $L(X)$ because $p$ is full. We can therefore deduce from the Rieffel correspondence that $I=C$. Thus $pZ$ is an $A$--$C$ imprimitivity bimodule. Similarly, $qZ$ is a $B$--$C$ imprimitivity bimodule. 

Note that the map $(x,qz)\mapsto x\cdot qz$ is bilinear, so there is a
well-defined map $\Omega$ on the algebraic tensor product $X\odot qZ$
satisfying \eqref{keyiso}, and which is $C$-linear. To see that it is $A$-linear, recall that the action of $A$ on $X$ is given by the product of the embedded copies in $L(X)$; thus for $a\in A$ and $x\in X$, we have
\[
\Omega(a\cdot x\otimes qz)=(a\cdot x)\cdot qz=a\cdot(x\cdot qz)=a\cdot\Omega(x\otimes qz).
\]
In the same way, the inner product $\langle y\,,\, x\rangle_B$ is given by the product $y^*x$ in $L(X)$, so
\begin{align*}
\langle x\otimes_B qz\,,\, y\otimes_B qw\rangle_C
&=\langle\langle y\,,\, x\rangle_B \cdot qz\,,\, qw\rangle_C
=\langle
(y^*x)\cdot qz\,,\, qw\rangle_C\\
&=\langle
x\cdot qz
\,,\, 
y\cdot qw\rangle_C
=\langle \Omega(x\otimes_B qz)\,,\, \Omega(y\otimes_B qw)\rangle_C,
\end{align*}
and $\Omega$ extends to an isometry of $(X\otimes_B qZ)_C$ into $(pZ)_C$.
To see that $\Omega$ has dense range and is therefore onto, note that
$L(X)$ acts nondegenerately on $Z$, so that
\begin{align*}
\range\Omega&\supset pL(X)q\cdot qZ= pL(X)q\cdot qL(X)\cdot Z\\
&=pL(X)qL(X)\cdot Z=pL(X)\cdot Z \\
&=pZ
\end{align*}
because $q$ is full.
Since $\Omega$ is a bimodule isomorphism which preserves the $C$-valued inner product, it must preserve the $A$-valued inner product as well.
\end{proof}

 To prove Theorem~\ref{thm-decomposition}, we apply Lemma~\ref{lem-ibm-isom} to the Combes bimodule $X\rtimes_{\gamma, r}G$ and a bimodule $Z$ coming from an application of Theorem~\ref{thm-sumup}. As it arises, $Z$ will be a left module over $L(X)\rtimes_r G$ rather than $L(X\rtimes_r G)$. Thus we shall have to identify $L(X)\rtimes_r G$ with $L(X\rtimes_r G)$. This is in some sense known, since Combes \emph{defined} $X\rtimes G$ to be the top right-hand corner in $L(X)\rtimes G$ \cite{com}; however, since we need to be very explicit about the identifications involved, we review the details.

We begin by recalling the construction of the Combes bimodule. Suppose as in the Theorem that ${}_{(A,\alpha)}(X,\gamma)_{(B,\beta)}$ is a Morita equivalence. For $f\in C_c(G,A), g\in C_c(G,B)$ and $z,w\in C_c(G, X)$ define
\begin{align}
f\cdot z(s)&=\int_G f(r)\cdot\gamma_r(z(r^{-1}s))\, dr\label{combesla}\\
z\cdot g(s)&=\int_G z(r)\cdot\beta_r(g(r^{-1}s))\, dr\label{combesra}\\
{}_{A\rtimes_\alpha G}\langle z\,,\, w\rangle(s)&=\int_G {}_A\langle z(r)\,,\, \gamma_s(w(s^{-1}r))\rangle\Delta_G(s^{-1}r)\, dr\label{combeslip}\\
\langle z\,,\, w\rangle_{B\rtimes_\beta G}(s)&=\int_G \beta_r^{-1}\big( \langle z(r)\,,\, w(rs)\rangle_B \big)\, dr.\label{combesrip}
\end{align}

\begin{prop}[Combes]\label{combesfull} With the above actions and inner products, $C_c(G,X)$ is a $C_c(G,A)$--$C_c(G,B)$ pre-imprimitivity bimodule whose completion is an $A\rtimes_\alpha G$--$B\rtimes_\beta G$ imprimitivity bimodule $X\rtimes_\gamma G$.
\end{prop}

\begin{proof}
We know from \cite[Proposition~3.5]{EKQR} that $C_c(G,X)$ is a pre-inner product module over $C_c(G,B)$, so we can complete to get a Hilbert $(B\rtimes_\beta G)$-module. We could also complete using the left-hand inner product, but the two satisfy the relation
\begin{equation}\label{ibidentity}
{}_{A\rtimes_\alpha G}\langle z\,,\, w\rangle \cdot y= 
z\cdot \langle w\,,\, y\rangle_{B\rtimes_\beta G},
\end{equation}
so the usual argument shows that the two semi-norms are equal, and the completions are the same. It follows that the left and right actions extend to actions of the full crossed products, and the relation \eqref{ibidentity} implies that the completion is an imprimitivity bimodule. 
\end{proof}

Composing functions with the identifications \eqref{idsinlinking} gives  embeddings 
$\iota_{11}$, $\iota_{12}$ and $\iota_{22}$ of $C_c(G,A)$, $C_c(G,X)$ and $C_c(G,B)$, respectively, in $C_c(G,L(X))$. The actions of $G$ on the corners combine to give an action $u$ of $G$ on $L(X)$, and then for the usual product in $L(X)\rtimes_u G$ we have
\[
\iota_{12}(z)^*\iota_{12}(w)=\iota_{22}(\langle z\,,\,w\rangle_{B\rtimes G})
\text{\ and\ }
\iota_{12}(z)\iota_{12}(w)^*=\iota_{11}({}_{A\rtimes G}\langle z\,,\, w\rangle)
\]
for $z,w\in C_c(G,X)$, and
\[
\iota_{11}(f)\iota_{12}(z)=\iota_{12}(f\cdot z)\text{\ and\ }
\iota_{12}(z)\iota_{22}(g)=\iota_{12}(z\cdot g)
\]
for $f\in C_c(G,A)$ and $g\in C_c(G, B)$. The Rieffel induction process $(X\oplus B)\dashind$ sets up a one-to-one correspondence $(\pi,U)\mapsto (L(\pi),L(U))$ between the covariant representations of $(B,G,\beta)$ and the covariant representations of $(L(X), G, u)$, and similarly for $(A,G,\alpha)$. Thus the embedding $\iota_{12}$ is isometric for the universal norm on $C_c(G,L(X))$, and hence extends to an embedding of $X\rtimes_\gamma G$ into $L(X)\rtimes_u G$. The image of $X\rtimes_\gamma G$ lies in the corner $\hat p\big(L(X)\rtimes_u G\big)\hat q$ associated to the images $\hat p$ and $\hat q$ of $p$ and $q$ under  the canonical embedding of $L(X)$ in $M(L(X)\rtimes_u G)$; indeed, since this range contains $\hat pC_c(G,L(X))\hat q$, it is precisely $\hat p\big(L(X)\rtimes_u G\big)\hat q$. Similar considerations show that $\iota_{11}$ and $\iota_{22}$ are isomorphisms onto the diagonal corners in $L(X)\rtimes_u G$. Thus with $\iota_{21}:=\widetilde\iota_{12}$, we have:

\begin{prop}\label{fulliso}
The maps $\iota_{ij}$ combine to give an isomorphism of 
\[
L(X\rtimes_\gamma G)=\begin{pmatrix} A\rtimes_\alpha G&X\rtimes_\gamma G\\(X\rtimes_\gamma G) \widetilde{\;}&B\rtimes_\beta G\end{pmatrix}
\]
onto $L(X)\rtimes_u G$, which carries $p_{L(X\rtimes_\gamma G)}$ and $q_{L(X\rtimes_\gamma G)}$ into $\hat p$ and $\hat q$.
\end{prop}

From this point of view, it is easy to deduce the analogues of Propositions~\ref{combesfull} and~\ref{fulliso} for the reduced crossed products. Let $\pi$ be a faithful representation of $B$ on $\H$. This induces a representation $L(\pi):=(X\oplus B)\dashind$ of $L(X)$ on $(X\otimes_B \H)\oplus \H$, and hence a representation $L(\pi)\widetilde{\;}\rtimes\lambda$ of $L(X)\rtimes_u G$ on $L^2(G, (X\otimes_B\H)\oplus\H)$. Moving this over to a representation of $L(X\rtimes_\gamma G)\cong L(X)\rtimes_u G$ on $L^2(G,X\otimes_B\H)\oplus L^2(G,\H)$ gives a representation $\rho$ of $L(X\rtimes_\gamma G)$ which restricts on the corners of $L(X\rtimes_\gamma G)$ to 
\[
\rho|_{A\rtimes_\alpha G}=\big((X\dashind\pi)\widetilde{\;}\rtimes\lambda\big)\oplus 0 \quad\text{and}\quad\rho|_{B\rtimes_\beta G}=0\oplus\big(\tilde\pi\rtimes\lambda\big).
\]
Thus $X\rtimes_{\gamma,r} G:=X\rtimes_\gamma G/\ker(\rho|_{X\rtimes G})$ implements a Morita equivalence between 
\[
A\rtimes_{\alpha,r} G= (A\rtimes_\alpha G)/(\ker (X\dashind\pi)\widetilde{\;}\rtimes\lambda) \mbox{\ and\ } B\rtimes_{\beta,r} G= (B\rtimes_\beta G)/\ker(\tilde\pi\rtimes\lambda).
\]
 Note that the representation $\rho$ induces an isomorphism of $L(X\rtimes_\gamma G)$ onto $\rho(L(X)\rtimes_u G)\cong L(X)\rtimes_{u,r} G$. We can view $B\rtimes_{\beta,r} G$ as the completion of $C_c(G,B)$ in the reduced norm $\|\cdot\|_r$, and $X\rtimes_{\gamma,r} G$ as the completion of $C_c(G,X)$ in the semi-norm defined by
\[
\|z\|^2:=\big\|\langle z\,,\,z\rangle_{B\rtimes_\beta G}\big\|_r.
\]
Thus:

\begin{cor}\label{rediso}
The maps $\iota_{ij}$ induce an isomorphism of $L(X\rtimes_{\gamma,r} G)$ onto $L(X)\rtimes_{u,r} G$.
\end{cor}

We now return to the situation of Theorem~\ref{thm-decomposition}. Recall that we seek an $L(X\rtimes_{\gamma,r} G)$--$D$ imprimitivity bimodule $Z$ to which we can apply Lemma~\ref{lem-ibm-isom}. We intend to find $Z$ by applying Theorem~\ref{thm-sumup} to the Morita equivalence ${}_{(L(X),u)}(X\oplus B)_{(B,\beta)}$ and identifying the left-hand algebra $L(X)\rtimes_{u,r} G$ with $L(X\rtimes_{\gamma,r}G)$ using Corollary~\ref{rediso}. Of course we have some checking to do:

\begin{lemma}\label{xplusb-pr}
Suppose that ${}_{(A,\alpha)}(X,\gamma)_{(B,\beta)}$ is proper with respect to the pre-imprimitivity bimodule ${}_{A_0}(X_0)_{B_0}$, and has generalized fixed-point algebra $D$. Let $B_1=\langle X_0\,,\, X_0\rangle_B$ and
\[
L(X_0):=\begin{pmatrix}A_0&X_0\\\widetilde{X_0}&B_1\end{pmatrix}.
\]
Then $\gamma\oplus \beta$ is proper with respect to ${}_{L(X_0)}(X_0\oplus B_1)_{B_1}$, and has generalized fixed-point algebra $D$. 
\end{lemma}

\begin{proof}Let  $x,y\in X_0$ and $b,c\in B_1$. For Definition~\ref{defn-pr}(1) we need to verify that 
\[
s\mapsto 
\!\!\!\!
{\phantom{\Big\rangle}}_{L(X)}
\Big\langle\begin{pmatrix}x\\b\end{pmatrix}\,,\,\begin{pmatrix}\gamma_s(y)\\\beta_s(c)\end{pmatrix}\Big\rangle
=\begin{pmatrix}
{}_A\langle x\,,\,\gamma_s(y)\rangle &x\cdot\beta_s(c^*)\\
\flat(\gamma_s(y)\cdot b^*)&b\beta_s(c^*)
\end{pmatrix}
\]
and its product with $\Delta(s)^{-1/2}$ are in $L^1(G,L(X))$. 
Since the action of $G$ is proper with respect to ${}_{A_0}(X_0)_{B_0}$ and ${}_{B_1}(B_1)_{B_1}$, it suffices to check that the functions
\begin{equation}\label{blah}
s\mapsto x\cdot\beta_s(c)\quad\text{and}\quad s\mapsto x\cdot\beta_s(c)\Delta(s)^{-1/2} 
\end{equation}
are integrable.
We know from Definition~\ref{defn-pr} for the  proper Morita equivalence  ${}_{(A_0,\alpha)}(X_0,\gamma)_{(B_0,\beta)}$
 that $s\mapsto{}_A\langle x\,,\,\gamma_s(y)\rangle$ and its product with $\Delta^{-1/2}$ are integrable. Thus the integrability of \eqref{blah} follows from the estimate
\[
\|x\cdot \beta_s(\langle y\,,\, z\rangle_B)\|\leq\|{}_A\langle x\,,\,\gamma_s(y)\rangle\cdot\gamma_s(z)\|\leq\| {}_A\langle x\,,\,\gamma_s(y)\rangle\|\,\|z\|.
\]

For Definition~\ref{defn-pr}(2), note that  $s\mapsto\gamma_s(x)\cdot b$ and $s\mapsto\beta_s(b)\cdot c$ and their products with $\Delta(s)^{-1/2}$ are integrable using Definition~\ref{defn-pr}(2) for ${}_{(A_0,\alpha)}(X_0,\gamma)_{(B_0,\beta)}$ and ${}_{(B_1,\beta)}(B_1,\beta)_{(B_1,\beta)}$.

To verify Definition~\ref{defn-pr}(3),  we  write $D'$ for the generalized fixed-point algebra associated to the action $\gamma\oplus\beta$, and $F:=(B,B_1)^\beta$, and define
\[
\Big\langle\begin{pmatrix}x\\b\end{pmatrix}\,,\,\begin{pmatrix}y\\c\end{pmatrix}\Big\rangle_{D'}:=\langle x\,,\,y\rangle_{D}+\langle b\,,\, c\rangle_{F}.
\]
Note that $D$ multiplies $B_1$,  and that the right-hand side belongs to $D$ because $F\subset D$ by part (1) of the Theorem.   Thus $D'=D$.  Straightforward calculations show that 
\[
\begin{pmatrix}z\\b'\end{pmatrix}\cdot \Big\langle\begin{pmatrix}x\\b\end{pmatrix}\,,\,\begin{pmatrix}y\\c\end{pmatrix}\Big\rangle_D\in \begin{pmatrix}X_0\\B_1\end{pmatrix},
\]
and that
\[
\int_G b'\beta_s\Big(\Big\langle\begin{pmatrix}x\\b\end{pmatrix}\,,\,\begin{pmatrix}y\\c\end{pmatrix}\Big\rangle_B\Big)\, ds=b'
\Big\langle\begin{pmatrix}x\\b\end{pmatrix}\,,\,\begin{pmatrix}y\\c\end{pmatrix}\Big\rangle_D
\]
 for $z\in X_0$ and $b'\in B_1$.
\end{proof}

At this point it is convenient to prove the assertion in Theorem~\ref{thm-decomposition} about saturated actions. 

\begin{prop} \label{lem-saturation} Suppose that the Morita equivalence ${}_{(A,\alpha)}(X,\gamma)_{(B,\beta)}$ is proper with respect to the pre-imprimitivity bimodule ${}_{A_0}(X_0)_{B_0}$. Then the following are equivalent:
\begin{enumerate}
\item[(1)] the action $\gamma$ is saturated with respect to $X_0$;
\item[(2)] the action $\beta$ is saturated with respect to $B_1$; and
\item[(3)] the action $\gamma\oplus\beta$ is saturated with respect to $X_0\oplus B_1$.
\end{enumerate}
\end{prop}

For the proof we need a standard lemma.

\begin{lemma}\label{lem-corner}
Let $E$ be an ideal in a $C^*$-algebra $C$ and let $p$ be a full projection in $M(C)$. Then $E=C$ if and only if $pEp=pCp$.
\end{lemma}

\begin{proof}
Recall that $pC$ is a $pCp$--$C$ imprimitivity bimodule. The result follows from the Rieffel correspondence: 
\begin{align*}
pC\dashind (E)&=\clsp\lbrace{}_{pCp}\langle pce\,,\, pd\rangle: c,d\in C,e\in E \rbrace\\
&=\clsp\lbrace pced^*p:c,d\in C,e\in E\rbrace\\
&=pEp
\end{align*}
is the corresponding ideal in $pCp$. 
\end{proof}

\begin{proof}[Proof of Proposition~\ref{lem-saturation}] 
As usual, we write  $E$ for the ideal in $L(X)\rtimes_{u,r} G$ spanned by functions of the form 
\begin{equation}\label{eq-la-fns}
s\mapsto \Delta(s)^{-1/2}\!\!\!\!
{\phantom{\Big\rangle}}_{L(X)}\Big\langle\begin{pmatrix}x\\b\end{pmatrix}\,,\,\begin{pmatrix}\gamma_s(y)\\\beta_s(c)\end{pmatrix}\Big\rangle
=\Delta(s)^{-1/2}
\begin{pmatrix}
{}_A\langle x\,,\,\gamma_s(y)\rangle &x\cdot\beta_s(c^*)\\
\flat(\gamma_s(y)\cdot b^*)&b\beta_s(c^*)
\end{pmatrix}
\end{equation}
where $x,y\in X_0$ and $b,c\in B_1$.
Since $L(X)\rtimes_{u,r} G\cong L(X\rtimes_{\gamma,r} G)$, two applications of Lemma~\ref{lem-corner} imply that
\[
E=L(X)\rtimes_{u,r} G \Longleftrightarrow \hat pL(X\rtimes_{\gamma,r}G)\hat p=A\rtimes_{\alpha,r} G \Longleftrightarrow \hat qL(X\rtimes_{\gamma,r}G)\hat q=B\rtimes_{\beta,r} G.
\]
Thus functions of the form~\eqref{eq-la-fns} are dense in $L(X)\rtimes_{u,r} G$ if and only if $A\rtimes_{\alpha,r} G$ is spanned by the functions $s\mapsto\Delta(s)^{-1/2}{}_A\langle x\,,\,\gamma_s(y)\rangle$ if and only if $B\rtimes_{\beta,r} G$ is spanned by the functions $s\mapsto\Delta(s)^{-1/2}b\beta_s(c^*)$. The result follows.
\end{proof}

\begin{proof}[Proof of Theorem~\ref{thm-decomposition}]
Parts (1) and (2) were proved earlier, and the statement about saturation is part of Proposition~\ref{lem-saturation}. We know from Lemma~\ref{xplusb-pr} and Proposition~\ref{lem-saturation} that $\gamma$ and $\gamma\oplus\beta$ are proper and saturated with respect to $X_0$ and ${}_{L(X_0)}(X_0\oplus B_1)_{B_1}$, respectively. Thus Theorem~\ref{thm-sumup}  gives two imprimitivity bimodules
\begin{equation}\label{eq-2modules}
{}_{B\rtimes_{\beta,r} G}(\overline{B_1})_{F}\quad\text{and}\quad{}_{L(X)\rtimes_{u,r} G}\overline{(X_0\oplus B_1)}_D,
\end{equation} 
where $F$ is an ideal of $D$. Since $L(X)\rtimes_{u,r} G\cong L(X\rtimes_{\gamma,r} G)$, we can apply Lemma~\ref{lem-ibm-isom} to the imprimitivity bimodules $X\rtimes_{\gamma,r} G$ and $\overline{X_0\oplus B_1}$. Thus to see the existence of the isomorphism, it suffices to prove that
\begin{equation}\label{eq-cornermodules}
\hat p(\overline{X_0\oplus B_1})\cong {}_{A\rtimes_{\alpha,r} G}\overline{(X_0)}_D \quad\text{\ and\ } \quad\hat q(\overline{X_0\oplus B_1})\cong{}_{B\rtimes_{\beta,r} G}(\overline{B_1})_D.
\end{equation}
 as imprimitivity bimodules; given this, it then follows from the  the Rieffel correspondence that $F=D$ because the imprimitivity bimodules in \eqref{eq-2modules}
 and \eqref{eq-cornermodules} based on $B_1$ are completed in the same norm.

That $\hat p(\overline{X_0\oplus B_1})$ and $\hat q(\overline{X_0\oplus B_1})$ are $(A\rtimes_{\alpha,r} G)$--$D$ and  $(B\rtimes_{\beta,r} G)$--$D$ imprimitivity bimodules, respectively, is proved in Lemma~\ref{lem-ibm-isom} (after again identifying $L(X)\rtimes_{u,r} G$ with $L(X\rtimes_{\gamma,r} G)$). 
Recall that  $E_0:={}_{L(X)\rtimes_{u,r} G}\langle X_0\oplus B_1\,,\,X_0\oplus B_1\rangle$ and that $E_0\cdot (X_0\oplus B_1)$ is dense in $(\overline{X_0\oplus B_1})_D$.
Since
\[
\hat p\!\!\!\!{\phantom{\Big\rangle}}_{L(X)\rtimes_{u,r} G}\bigg\langle\begin{pmatrix}x\\b\end{pmatrix}\,,\,\begin{pmatrix}y\\c\end{pmatrix}\bigg\rangle=
\!\!\!\!{\phantom{\Big\rangle}}_{L(X)\rtimes_{u,r} G}\bigg\langle\begin{pmatrix}x\\0\end{pmatrix}\,,\,\begin{pmatrix}y\\c\end{pmatrix}\bigg\rangle
\]
implies 
 \[
\hat p\!\!\!\!{\phantom{\Big\rangle}}_{L(X)\rtimes_{u,r} G}\bigg\langle\begin{pmatrix}x\\b\end{pmatrix}\,,\,\begin{pmatrix}y\\c\end{pmatrix}\bigg\rangle
\cdot \begin{pmatrix}z\\d\end{pmatrix}=
\begin{pmatrix}x\\0\end{pmatrix}\cdot\bigg\langle\begin{pmatrix}y\\c\end{pmatrix}\,,\,\begin{pmatrix}z\\d\end{pmatrix}\bigg\rangle_D\in\begin{pmatrix}X_0\\\lbrace 0\rbrace\end{pmatrix},
\]
we obtain that $\hat p(\overline{X_0\oplus B_1})_D=(\overline{X_0\oplus \lbrace 0\rbrace})_D$. That $\hat p(\overline{X_0\oplus B_1})\cong {}_{A\rtimes_{\alpha,r} G}\overline{(X_0)}_D$ is now clear because the inclusion of $X_0$ into $(\overline{X_0\oplus \lbrace 0\rbrace})_D$ preserves both inner products and the $D$-action. Similarly, $\hat q\overline{(X_0\oplus B_1)}\cong{}_{B\rtimes_{\beta,r} G}(\overline{B_1})_D$.

Finally, to get the formula for the isomorphism, we need to chase through our identifications.
Here, $f\cdot b$ means the left action of $f\in X\rtimes_{\gamma,r} G\subset L(X\rtimes_{\gamma,r} G)\cong L(X)\rtimes_{u,r} G$ on $b\in B_1\subset X_0\oplus B_1$.
 Thus we have a formula for the action provided 
\[
\begin{pmatrix}0&f\\0&0\end{pmatrix}=\!\!\!\!{\phantom{\Big\rangle}}_{L(X)\rtimes_{u,r} G}\Big\langle\begin{pmatrix}x\\0\end{pmatrix}\,,\,\begin{pmatrix}0\\c\end{pmatrix}\Big\rangle,
\] 
which means that $f$ must have the form $f(s)=x\cdot\beta_s(c^*)\Delta(s)^{-1/2}$. If so,
\begin{align*}
\begin{pmatrix}f\cdot b\\0\end{pmatrix}=\begin{pmatrix}0&f\\0&0\end{pmatrix}\cdot\begin{pmatrix}0\\b\end{pmatrix}
&=\int_G \begin{pmatrix}0&f\\0&0\end{pmatrix}(s)\cdot\gamma_s\oplus\beta_s\begin{pmatrix}0\\b\end{pmatrix}\Delta(s)^{1/2}\,ds\\
&=\int_G
\begin{pmatrix}0&f(s)\\0&0\end{pmatrix}\begin{pmatrix}0\\\beta_s(b)\end{pmatrix}\Delta(s)^{1/2}\,ds\\
&=\int_G\begin{pmatrix}f(s)\cdot\beta_s(b)\Delta(s)^{1/2}\\0\end{pmatrix}\,ds\\
&=\begin{pmatrix}\int_G f(s)\cdot\beta_s(b)\Delta(s)^{1/2}\,ds\\0\end{pmatrix},
\end{align*}
which gives the right formula.
\end{proof}

\section{The symmetric imprimitivity theorem for induced algebras}\label{sec-sit}

We recall the set-up of the symmetric imprimitivity theorem.  We start 
with commuting free and proper actions of $K$ and $H$ on the left and right of 
a locally compact space $P$, 
and commuting actions $\alpha:K\to\Aut A$ and $\beta:H\to\Aut A$ on the same 
$C^*$-algebra $A$. The actions of $K$ and $H$ induce actions $\lt$ and $\rt$ on $C_0(P)$, so that $\lt_t(f)(p)=f(t^{-1}\cdot p)$ and $\rt_s(f)(p)=f(p\cdot s)$.

The \emph{induced $C^*$-algebra} $\Ind \alpha=\Ind_K^P(A,\alpha)$ is the 
$C^*$-subalgebra of $C_b(P,A)$ consisting of the functions $f$ such that 
$f(t\cdot p)=\alpha_t(f(p))$ for $t\in K$ and $p\in P$, and such that the 
function $K\cdot p\mapsto\|f(p)\|$ belongs to $C_0(K\backslash P)$.  The diagonal 
action $\tau:=\rt\otimes\beta$ of $H$ on $C_b(P,A)\subset M(C_0(P,A))$ restricts 
to a strongly continuous 
action of $H$ on $\Ind\alpha$ which is characterized by
\[
\tau_s(f)(p)=(\rt\otimes\beta)_s(f)(p)=\beta_s(f(p\cdot s))
\]
(see \cite[Lemma~5.1]{hrw}).
Likewise, $\Ind\beta$ consists of the bounded continuous functions $f:P\to A$ 
such that $f(p\cdot s)=\beta_s^{-1}(f(p))$ for $s\in H$ and $p\in P$, and 
such that $p\cdot H\mapsto \|f(p)\|$ vanishes at infinity on $P/H$, and there is
a natural action $\sigma:=\lt\otimes\alpha$ of $K$ on $\Ind\beta$ given by
\[
\sigma_t(f)(p)=(\lt\otimes\alpha)_t(f)(p)=\alpha_t(f(t^{-1}\cdot p)).
\]

The symmetric imprimitivity theorem of \cite{rae} says that $Z_0:=C_c(P,A)$ completes 
to an $(\Ind\beta\rtimes_\sigma K)$--$(\Ind\alpha\rtimes_\tau H)$ imprimitivity 
bimodule $Z$, in which the actions and inner products are given by\footnote{These are the 
same as \cite[Equations~(4.14--4.17)]{hrw}, modulo replacing $h\in H$ with 
$s\in H$, $s\in G$ 
with $t\in K$, $x\in X$ with $p\in P$, and correcting the mistake in 
\cite[Equation~(4.15)]{hrw} by replacing $c(h^{-1}, x\cdot h^{-1})$ by 
$c(h, x\cdot h^{-1})$.}
\begin{align}
b\cdot f(p)&=\int_K b(t,p)\alpha_t\big( f(t^{-1}\cdot p)\big)\Delta_K(t)^{1/2}\, 
dt\label{sit-la}\\
f\cdot c(p)&=\int_H \beta_s^{-1}\big( f(p\cdot s^{-1})c(s,p\cdot s^{-1})  \big)
\Delta_H(s)^{-1/2}\, ds\label{sit-ra}\\
{}_{(\Ind\beta)\rtimes K}\langle f\,,\, g\rangle(t,p)&=\Delta_K(t)^{-1/2}
\int_H \beta_s\big( f(p\cdot s)\alpha_t(g(t^{-1}\cdot p\cdot s)^*) \big)\, ds
\label{sit-lip}\\
\langle f\,,\, g\rangle_{(\Ind\alpha)\rtimes H}(s,p)&=\Delta_H(s)^{-1/2}
\int_K\alpha_t\big(f(t^{-1}\cdot p)^*\beta_s(g(t^{-1}\cdot p\cdot s))  \big)
\, dt\label{sit-rip},
\end{align}
for $f,g\in C_c(P,A)$, $b\in C_c(K, \Ind\beta)$ and $c\in C_c(H, \Ind\alpha)$.
Quigg and Spielberg proved in \cite{qs} that $Z_0$, with the same 
formulas~\eqref{sit-la}--\eqref{sit-rip}, also completes to give a
Morita equivalence $Z_r$ between the reduced crossed products $(\Ind\beta)\rtimes_{\sigma,r} K$ and 
$(\Ind\alpha)\rtimes_{\tau,r} H$ (see also 
\cite[Corollary~2]{aHRqs}). 

The Morita equivalence of $(\Ind\beta)\rtimes_{\sigma,r} K$ and 
$(\Ind\alpha)\rtimes_{\tau,r} H$ has been derived in different fashion by 
Kasparov \cite[Theorem~3.15]{kas}: his equivalence is implemented by a bimodule
$(X\rtimes_r H)\otimes_{C_0(P,A)\rtimes_r(H\times K)}(Y\rtimes_r K)$ in which 
$X$ and $Y$ are one-sided analogues of the imprimitivity bimodule $Z$. 
Here we shall show that the machine of \S\ref{genRieffel} and \S\ref{decomp} 
gives an isomorphism between the bimodules of Quigg-Spielberg and Kasparov.

We intend to apply Theorem~\ref{thm-decomposition} to the 
$(\Ind\beta)$--$(C_0(P,A)\rtimes_{\tau} H)$ imprimitivity bimodule $X$ 
obtained by ignoring the 
action of $K$ in \eqref{sit-la}--\eqref{sit-rip}. (This is isomorphic to the 
dual of the bimodule constructed in \cite{rw85}; see Remark~\ref{cfaHRW}.) In this situation we have 
\[
C_0(P,A)\rtimes_\tau H=C_0(P,A)\rtimes_{\tau,r}H
\]
 by \cite[Corollary~3]{aHRqs}, 
so we do not have to take a quotient of $X$ to handle reduced crossed products. 
For the pre-imprimitivity 
bimodule ${}_{A_0}(X_0)_{B_0}$, we take the bimodule based on $X_0:=C_c(P,A)$. We
therefore have to verify that the diagonal action $\sigma=\lt\otimes\alpha$ on $X$ 
implements a Morita equivalence and is proper with respect to $X_0$. To avoid clumsy
notation, we write $\sigma$ for the various actions on $X_0$, $C_0(P,A)$ and 
$\Ind\beta$ induced by $\lt\otimes\alpha$. We also write 
$C$ for $C_0(P,A)$. 

\begin{lemma}\label{checkproper} 
The diagonal action $\sigma=\lt\otimes\alpha$ induces continuous actions $\sigma$ 
of $K$ on 
$\Ind\beta$ and on $X$ and $\sigma \rtimes\id$ on $C\rtimes_{\tau}H$ such that
$(X,K,\sigma)$ is a Morita equivalence between $(\Ind\beta, K,\sigma)$ and 
$(C\rtimes_{\tau}H,K,\sigma\rtimes\id)$. The action of $K$ on $X$ is proper 
and saturated with respect to the pre-imprimitivity bimodule 
$X_0={}_{\Ind\beta}C_c(P,A)_{C_c(H,C)}$.
\end{lemma}

\begin{proof}
A compactness argument shows that the action $\sigma\rtimes\id$ is continuous on 
$C_c(H,C)$ for the inductive limit topology, hence on $C\rtimes H$. 
The action on $\Ind\beta$ is continuous by 
\cite[Lemma~5.1]{hrw}. To see that $\sigma$ is continuous on $X$, note first 
that for any $y\in X_0$, \[
\|y\|^2=\sup_{p\in \supp y}\big\|{}_{\Ind\beta}\langle y\,,\,y\rangle(p)\big\|
\leq \sup_{p\in \supp y}\int_H \|y(p\cdot s)y(p\cdot s)^*\|\,ds.
\]
using Equation~\eqref{sit-lip}.
Because $H$ acts properly, $L:=\{s:(\supp y)\cdot s\cap\supp y\not=\emptyset\}$ 
is compact, and we have $\|y\|^2\leq \mu(L)\|y\|_\infty^2$. Now let $x\in X_0$, use 
a compactness argument to see that $\sigma_t(x)-x\to 0$  uniformly with 
support in a fixed compact neighbourhood as $t\to e$, and take $y=\sigma_t(x)-x$ 
in the inequality to see that $\|\sigma_t(x)-x\|\to 0$ as $t\to e$. Thus $\sigma$ 
is continuous on $X$. It is easy to check that the triple 
$(\sigma,\sigma,\sigma\rtimes\id)$ has the required algebraic properties, 
and hence $(X,K,\sigma)$ is a Morita equivalence, as claimed.

Because
\begin{align*}
{}_{\Ind\beta}\langle x\,,\,\sigma_t(y)\rangle (p)
&=\int_H \beta_s(x(p\cdot s)\sigma_t(y)(p\cdot s)\, ds\\
&=\int_H \beta_s\big(x(p\cdot s)\alpha_t(y(t^{-1}\cdot p\cdot s))\big)\,ds
\end{align*}
vanishes unless $(\supp x)\cap t\cdot(\supp y)$ is nonempty, the function 
$t\mapsto {}_{\Ind\beta}\langle x\,,\,\sigma_t(y)\rangle$ has compact support. 
It is continuous because $\sigma$ is, and hence the integrability conditions 
in Definition~\ref{defn-pr}(1) are trivially satisfied. Similar considerations 
show that Definition~\ref{defn-pr}(2) holds.

To verify the existence of the multiplier $\langle x\,,\,y\rangle_D$ of $C\rtimes_\tau H$, we consider 
the function $\langle x\,,\,y\rangle_D$ defined by the right-hand side of 
\eqref{sit-rip}. This is a continuous function of compact support from $H$ to $\Ind\alpha$.
 Since the inclusion of $\Ind\alpha$ in $M(C)=M(C_0(P,A))$ is equivariant 
for the actions $\tau$, it induces a homomorphism of $(\Ind\alpha)\rtimes_\tau H$ 
into $M(C\rtimes_\tau H)$; we take $\langle x\,,\,y\rangle_D$ to be the multiplier 
defined by the function in $C_c(H,\Ind \alpha)\subset (\Ind\alpha)\rtimes_\tau H$. 
For $b\in C_c(H,C)$, the product $b\langle x\,,\,y\rangle_D$ 
is given by the usual convolution formula
\[
(b\langle x\,,\,y\rangle_D)(s)=\int_H b(r)\tau_r(\langle x\,,\,y\rangle_D)(r^{-1}s)\,dr,
\]
which has compact support because $b$ and $\langle x\,,\,y\rangle_D$ do. 
Because evaluation at $p\in P$ is a homomorphism on $C$, it pulls through the integral, and
\begin{equation}\label{verifypropip}
(b\langle x\,,\,y\rangle_D)(s,p)=
\int_H b(r,p)\beta_r(\langle x\,,\,y\rangle_D(r^{-1}s,p\cdot r))\,dr.
\end{equation}
But now we notice that
\begin{align*}
\langle x\,,\,y\rangle_D(s,p)
&=\int_K \alpha_t\big(\langle x\,,\,y\rangle_{C_c(H,C)}(s,t^{-1}\cdot p)\big)\,dt\\
&=\int_K(\sigma\rtimes\id)_t(\langle x\,,\,y\rangle_{C_c(H,C)})(s,p)\,dt,
\end{align*}
and changing the order of integration in \eqref{verifypropip} shows that
\[
(b\langle x\,,\,y\rangle_D)(s,p)=
\int_K \big(b(\sigma\rtimes\id)_t(\langle x\,,\,y\rangle_{C_c(H,C)})\big)(s,p)\,dt.
\]
Equation \eqref{pr-eq} now follows from Lemma~\ref{pullpin}.

Pulling the variable $p$ through the integrals (see Remark~\ref{pullpinbimods}) 
shows that the action of
$C\rtimes_\tau H$ on $X$ is the integrated form of the action $M$ of $C=C_0(P,A)$ 
by pointwise multiplication on $X_0$ and the action of $H$ defined by 
$(z\cdot s)(p)=\beta_s^{-1}(z(p\cdot s^{-1}))$. Since the restriction of a
representation $\pi\rtimes U$ of $C\rtimes_\tau H$ to $(\Ind\alpha)\rtimes_\tau H$ is 
$\bar \pi|_{\Ind\alpha}\rtimes U$ and $\overline{M}|_{\Ind \alpha}$ is again 
pointwise multiplication, it follows from Lemma~\ref{pullpin} that the action
of $(\Ind\alpha)\rtimes_\tau H$ is given by the same formula \eqref{sit-ra}. 
This implies that $z\cdot\langle x\,,\,y\rangle_D$ also has compact support in $P$, 
and hence belongs to $X_0=C_c(P,A)$. This completes the proof that $\sigma$ 
acts properly.

To verify that $\sigma$ is saturated, we note that
\begin{align}\label{lipsagree}
{}_E\langle x\,,\,y\rangle(t,p)
&=\Delta_K(t)^{-1/2}{}_{\Ind\beta}\langle x\,,\,\sigma_t(y)\rangle(p)\\
&=\Delta_K(t)^{-1/2}\int_H \beta_s\big(x(p\cdot s)\sigma_t(y)
(p\cdot s)^*\big)\,ds\notag\\
&={}_{(\Ind\beta)\rtimes K}\langle x\,,\,y\rangle(t,p).\notag
\end{align}
Thus $E_0$ is the subalgebra of $C_c(K,\Ind \beta)$ spanned by the range 
of the inner product \eqref{sit-lip}; since the symmetric imprimitivity theorem 
asserts, \emph{inter alia}, that this is dense in $(\Ind \beta)\rtimes_\sigma K$, 
we deduce that $\sigma$ is saturated with respect to $X_0$.
\end{proof}

Lemma~\ref{checkproper} allows us to apply Theorem~\ref{thm-decomposition}
to the action $\sigma$ of $K$ on ${}_{\Ind\beta}X_{C\rtimes H}$. This gives 
a $(C\rtimes_r(H\times K))$--$(C\rtimes_r H)^{\sigma\times\id}$ 
imprimitivity bimodule $\overline{B_1}$, an $((\Ind\beta)\rtimes_r K)$--
$(C\rtimes_r H)^{\sigma\times\id}$ imprimitivity bimodule $\overline{X_0}$, 
and an isomorphism
\[
(X\rtimes_{\sigma,r} K)\otimes_{C\rtimes_r(H\times K)} \overline{B_1}\cong 
\overline{X_0}.
\]
The Combes bimodule $X\rtimes_{\sigma,r} K$ is one of the bimodules 
appearing in Kasparov's tensor product equivalence. We will
identify  $\overline{X_0}$ with the bimodule 
$Z_r$ of the reduced symmetric imprimitivity theorem.
To identify $\overline{B_1}$ with a familiar object, we need the 
$(C\rtimes_{\sigma,r} K)$--$(\Ind \alpha)$ imprimitivity bimodule $Y$ 
obtained by ignoring the 
action of $H$ in \eqref{sit-la}--\eqref{sit-rip}; this is 
the bimodule constructed in \cite{rw85}. The mirror image of 
Lemma~\ref{checkproper} says that the action $\tau$ of $H$ on 
$Y_0:=C_c(P,A)$ induces a Morita equivalence $(Y,H,\tau)$ between
$(C\rtimes_{\sigma,r} K,H,\tau\rtimes\id)$ and $(\Ind\alpha, H,\tau)$. We will 
show that $\overline{B_1}=Y\rtimes_{\tau,r} H$, and deduce:

\begin{thm}\label{astridsiso} 
There is an isomorphism
\[
\Omega:(X\rtimes_{\sigma,r}K)\otimes_{C\rtimes_r(H\times K)}
(Y\rtimes_{\tau,r}H)\to Z_r
\]
of $(\Ind\beta\rtimes_{\sigma,r}K)$--
$(\Ind\alpha\rtimes_{\tau,r}H)$ imprimitivity bimodules such that
$\Omega(f\otimes b)(p)$ is given by
\begin{equation}
=\int_K\int_H\beta_s^{-1}\big( f(t,p\cdot s^{-1})\big) 
\beta_s^{-1}\alpha_t\big(b(s,t^{-1}\cdot p\cdot s^{-1})  \big)
\Delta_K(t)^{1/2}\Delta_H(s)^{-1/2}\, ds\, dt.\label{eq-must-resolve}
\end{equation}
for $f\in C_c(K\times P,A))$ and $b\in C_c(H\times P,A)$.
\end{thm}

The next lemma will help us identify $X_0$ and $\overline{B_1}$.

\begin{lemma}\label{idibs}
Suppose that ${}_{C_0}(Y_0)_{D_0}$ and
${}_{A_0}(X_0)_{B_0}$ are pre-imprimitivity bimodules, and that we 
have $*$-homomorphisms $\phi:C_0\to A_0$ and $\psi:D_0\to B_0$ which 
extend to isomorphisms of $C=\overline{C_0}$ onto $A=\overline{A_0}$ 
and $D$ onto $B$.  
If $\iota:Y_0\to X_0$ is a linear map such that either
\begin{enumerate}
\item[(1)] $\iota(c\cdot y)=\phi(c)\cdot\iota(y)$, 
$\iota(y\cdot d)=\iota(y)\cdot\psi(d)$ and 
$\psi(\langle y\,,\,z\rangle_{D_0})=\langle \iota(y)\,,\,\iota(z)\rangle_{B_0}$, or

\item[(2)] $\iota(y\cdot d)=\iota(y)\cdot\psi(d)$,
$\psi(\langle y\,,\,z\rangle_{D_0})=\langle \iota(y)\,,\,\iota(z)\rangle_{B_0}$ 
and $\phi({}_{C_0}\langle y\,,\,z\rangle)={}_{A_0}\langle \iota(y)\,,\,\iota(z)\rangle$,
\end{enumerate}
then $(\phi,\iota,\psi)$ extends to an isomorphism of the imprimitivity bimodule 
${}_C\overline{Y_0}_D$ onto ${}_A\overline{X_0}_B$.
\end{lemma}

\begin{proof}
Suppose (1) is satisfied. Then
\[
\|\iota(y)\|^2=\|\langle \iota(y)\,,\,\iota(y)\rangle_{B_0}\|
=\|\psi(\langle y\,,\,y\rangle_{D_0})\|=\|\langle y\,,\,y\rangle_{D_0}\|=\|y\|^2,
\]
so $\iota$ extends to an isometry of $\overline{Y_0}$ onto a closed subspace 
$\iota(\overline{Y_0})$ of $\overline{X_0}$; the properties in (1) imply
that $\iota(\overline{Y_0})$ is an $A$--$B$ submodule of $\overline{X_0}$. Since
\begin{align*}
\phi({}_{C_0}\langle y\,,\,z\rangle)\cdot \iota(w)
&=\iota({}_{C_0}\langle y\,,\,z\rangle\cdot w)=\iota(y\cdot \langle z\,,\,w\rangle_{D_0})\\
&=\iota(y)\cdot \psi(\langle z\,,\,w\rangle_{D_0})
=\iota(y)\cdot \langle \iota(z)\,,\, \iota(w)\rangle_{B_0}\\
&={}_{A_0}\langle \iota(y)\,,\,\iota(z)\rangle\cdot\iota(w),
\end{align*}
the triple $(\phi,\iota,\psi)$ extends to an imprimitivity-bimodule homomorphism. 
Since the ranges of the inner products on $Y_0$ span dense ideals in $C$ and $D$,
and since $\phi$ and $\psi$ are isomorphisms on $C$ and $D$, the ranges of 
the inner products on $\iota(\overline{Y_0})$ span dense ideals in $A$ and $B$. 
Thus it follows from the Rieffel correspondence that 
$\iota(\overline{Y_0})=\overline{X_0}$.

If (2) holds, the map $\iota$ extends as before, but now we use that 
$(\phi,\iota,\psi)$ preserves the inner products to see that 
$\iota(c\cdot y)=\phi(c)\cdot\iota(y)$ for $c$ of the form ${}_{C_0}\langle w\,,\,z\rangle$, 
and extend this by continuity to $y\in\overline{Y_0}$ and $c\in C$.
Now (1) applies.
\end{proof}

\begin{remark}
If $\iota:Y_D\to X_B$ is a surjective linear map between Hilbert modules and $\psi:D\to B$ is a homomorphism, then
\[
\psi(\langle y\,,\,z\rangle_D)=\langle \iota(y)\,,\,\iota(z)\rangle_B\mbox{ for all $y,z\in Y$ }
\Longrightarrow \iota(y\cdot d)=\iota(y)\cdot\psi(d).
\]
Indeed, 
\begin{align*}
\langle\iota(z)\,,\,\iota(y\cdot d)\rangle_B
&=\psi(\langle z\,,\, y\cdot d\rangle_D)=\psi(\langle z\,,\, y\rangle_D\cdot d)\\
&=\psi(\langle z\,,\, y\rangle_D)\cdot \psi(d)=\langle \iota(z)\,,\,\iota(y)\rangle_B\cdot \psi(d)\\
&=\langle \iota(z)\,,\,\iota(y)\cdot \psi(d)\rangle_B,
\end{align*}
which implies $\iota(y\cdot d)=\iota(y)\cdot\psi(d)$ because every $x\in X$ has the form $\iota(z)$. 

It is tempting to believe that we can therefore remove one of the hypotheses in both (1) and (2) of Lemma~\ref{idibs}. However, it is important in our applications of Lemma~\ref{idibs} that we can deduce surjectivity of $\iota$ from it. If we already knew that the range of $\iota$ were dense in $X_B$, we would not need to check, for example, that $\iota(y\cdot d)=\iota(y)\cdot\psi(d)$ in (1).  
\end{remark}

\begin{proof}[Proof of Theorem~\ref{astridsiso}]
To identify $\overline{X_0}$ with $Z_r$, recall from the proof of 
Lemma~\ref{checkproper} that the multipliers $\langle x,y\rangle_D$ which span
$D_0\subset M(C\rtimes_{\tau}H)$ are the images of 
$\langle x,y\rangle_{(\Ind \alpha)\rtimes_\tau H}$ under the natural homomorphism
of $(\Ind \alpha)\rtimes_\tau H$ into $M(C\rtimes_{\tau}H)$ induced by the 
equivariant embedding of $\Ind\alpha$ in $M(C)=M(C_0(P,A))$. A faithful
nondegenerate representation of $C$ extends to a faithful nondegenerate
representation $\bar\pi$ of $\Ind \alpha$, and extending the regular 
representation $\widetilde\pi\times\lambda$ to $M(C\rtimes_\tau H)$ and 
restricting gives the regular representation $\widetilde{\bar\pi}\times\lambda$
of $(\Ind\alpha)\rtimes_\tau H$. Thus $(\Ind \alpha)\rtimes_{\tau,r} H$ embeds 
faithfully in $M(C\rtimes_\tau H)$, with range the generalized fixed-point algebra 
$(C\rtimes_\tau H)^{\sigma\rtimes\id}:=\overline{D_0}$. When we view 
$(C\rtimes_\tau H)^{\sigma\rtimes\id}$ as $(\Ind \alpha)\rtimes_{\tau,r} H$, the 
right-hand inner product on $X_0=C_c(P,A)$ becomes the inner product
\eqref{sit-rip} on $Z_0$, and we have already observed in the proof of
Lemma~\ref{checkproper} that the right action is given by \eqref{sit-ra}.
The calculation \eqref{lipsagree} shows that the left inner product agrees too, 
and it follows from Lemma~\ref{idibs} that $\overline{X_0}=Z_r$ as 
imprimitivity bimodules. 

The algebra $B_1$ is the dense subalgebra of $B:=C\rtimes_\tau H=C\rtimes_{\tau,r} H$ spanned by 
the functions of the form 
$\langle x\,,\,y\rangle_B:(s,p)\mapsto \Delta_H(s)^{-1/2}x(p)^*\beta_s(y(p\cdot s))$
for $x,y\in C_c(P,A)$. These functions belong to $C_c(H\times P,A)$;
we aim to use Lemma~\ref{idibs} to prove that the inclusion $\iota$ of 
$B_1$ in $C_c(H\times P,A)$
extends to an isomorphism of ${}_{B\rtimes_{\sigma\rtimes\id,r}K}(\overline{B_1})_D$ onto the Combes bimodule 
\[
{}_{(C\rtimes_\sigma K)\rtimes_{\tau\rtimes\id,r}H}(Y\rtimes_{\tau,r}H)_{(\Ind\alpha)\rtimes_{\tau,r} H}.
\]
For the map $\psi$ we use the identification of
$D=(C\rtimes_\tau H)^{\sigma\rtimes\id}$ with $(\Ind \alpha)\rtimes_{\tau, r} H$.
The action of $\langle b\,,\,c\rangle_D\in C_c(H,\Ind\alpha)$ on $C_c(H,C)$ is by
the usual formula for convolution in $C_c(H,C)$, so
\begin{equation*}
(a\cdot \langle b\,,\,c\rangle_D)(s)=
\int_H a(r)\tau_r\big(\langle b\,,\,c\rangle_D)(r^{-1}s)\big)\,dr \ \mbox{ for $a,b,c\in B_1$}, 
\end{equation*} 
which agrees with the action of $C_c(H,\Ind\alpha)$ on $Y\rtimes_{\tau,r} H$ 
given by \eqref{combesra}. In other words, $\iota(a\cdot d)=\iota(a)\cdot \psi(d)$ for $d\in D_0$.

To prove that $\psi(\langle b\,,\,c\rangle_{B^{\sigma\rtimes\id}})
=\langle \iota(b)\,,\,\iota(c)\rangle_{C_c(H,\Ind\alpha)}$, we show that
\[
\iota(a)\cdot \psi(\langle b\,,\,c\rangle_{B^{\sigma\rtimes\id}})=
\iota(a\cdot\langle b\,,\,c\rangle_{B^{\sigma\rtimes\id}})=
\iota(a)\cdot \langle \iota(b)\,,\,\iota(c)\rangle_{C_c(H,\Ind\alpha)}
\]
for every $a\in B_1$. From the characterizing property of the 
$B^{\sigma\rtimes\id}$-valued inner product associated to the proper action $\sigma\rtimes\id$ on $B$, and from Lemma~\ref{pullpin}, we have
\begin{align*}
\big(a\cdot\langle b\,,\,c\rangle_{B^{\sigma\rtimes\id}}\big)(s)
&=\int_K\big(a(\sigma\rtimes\id)_t(b^*c)\big)(s)\,dt\\
&={\int_K}\int_H a(r)\tau_r\sigma_t(b^*c(r^{-1}s))\,dr\,dt\\
&=\int_K\int_H\int_H a(r)\tau_r\sigma_t
\big(b^*(u)\tau_u(c(u^{-1}r^{-1}s))\big)\,du\,dr\,dt\\
&=\int_K\int_H\int_H a(r)\tau_{ru}\sigma_t
\big(b(u^{-1})^*c(u^{-1}r^{-1}s)\big)\Delta_H(u)^{-1}\,du\,dr\,dt.
\end{align*}
On the other hand, from \eqref{combesra} we have
\begin{align*}
\big(\iota(a)\cdot \langle \iota(b)\,,\,\iota(c)\rangle_{C_c(H,\Ind\alpha)}\big)(s)
&=\int_H a(r)\tau_r\big(\langle\iota(b)\,,\,\iota(c)\rangle_{C_c(H,\Ind\alpha)}(r^{-1}s)\big)\,dr\\
&=\int_H\int_H a(r)\tau_{rv^{-1}}
\big(\langle b(v)\,,\,c(vr^{-1}s)\rangle_{\Ind\alpha}\big)\,dv\,dr\\
&=\int_H\int_H\int_K a(r)\tau_{rv^{-1}}\sigma_t(b(v)^*c(vr^{-1}s))\,dt\,\,dv\,dr,
\end{align*}
which on substituting $u=v^{-1}$ gives 
$(a\cdot\langle b\,,\,c\rangle_{B^{\sigma\rtimes\id}})(s)$. 
Thus $(\iota,\psi)$ preserves the right inner product.

The left inner product on $B_1$ is given by
\begin{align*}
{}_{B\rtimes_r K}\langle b\,,\,c\rangle(t,s)
&=\Delta_K(t)^{-1/2}\big(b(\sigma\rtimes\id)_t(c^*)\big)(s)\\
&=\Delta_K(t)^{-1/2}\int_H b(r)\tau_r\sigma_t(c^*(r^{-1}s))\,dr\\
&=\Delta_K(t)^{-1/2}\int_H b(r)\sigma_t\tau_r\tau_{r^{-1}s}
(c(s^{-1}r)^*)\Delta_H(r^{-1}s)^{-1}\,dr\\
&=\int_H{}_{C\rtimes K}\langle b(r)\,,\,\tau_s(c(s^{-1}r))\rangle(t)
\Delta_H(s^{-1}r)\,dr,
\end{align*}
which is the left inner product of $b,c\in C_c(H,Y_0)$ described in 
\eqref{combeslip}. In other words, with 
$\phi:C_c(K\times H,C)\to C_c(H\times K,C)$ defined by 
$\phi(z)(t,s)=z(s,t)$, we have $\phi({}_{B\rtimes_r K}\langle b\,,\,c\rangle)
={}_{(C\rtimes K)\rtimes_r H}\langle \iota(b)\,,\,\iota(c)\rangle$.

Since the homomorphisms $\phi$, $\psi$ certainly extend to isomorphisms 
on the completions, we can deduce from Lemma~\ref{idibs} that $(\phi,\iota,\psi)$ 
extends to an isomorphism of $\overline{B_1}$ onto $Y\rtimes_{\tau,r}H$.

It remains to verify the formula for the isomorphism. Using 
\eqref{eq-isom} and then \eqref{sit-ra} gives 
\begin{align*}
\Omega(f\otimes b)&=\int_K f(t)\cdot (\sigma\rtimes\id)_t(b) \Delta_K(t)^{1/2}\,dt\\
&=\int_K\int_H \tau_s^{-1}\big( f(t)\cdot \sigma_t(b(s))  
\big)\Delta_K(t)^{1/2}\Delta_H(s)^{-1/2}\, ds\, dt
\end{align*}
when  $b\in B_1$ and $f\in L^1(K,X)$ has the form $f(t)=
\Delta_K(t)^{-1/2}x\cdot (\sigma\rtimes\id)_t(c^*)$ for $x\in X_0$ and $c\in B_1$.
Inserting the variable $p$ (see Remark~\ref{pullpinbimods}) gives
\begin{align*}
\Omega(f&\otimes b)(p)
=\int_K\int_H  \tau_s^{-1}(f(t))(p)\tau_s^{-1}\sigma_t(b(s))(p)
\Delta_K(t)^{1/2}\Delta_H(s)^{-1/2}\, ds\, dt\\
&=\int_K\int_H \beta_s^{-1}\big( f(t,p\cdot s^{-1}) \big)
\beta_s^{-1}\alpha_t\big( b(s,t^{-1}\cdot p\cdot s^{-1}) \big)
\Delta_K(t)^{1/2}\Delta_H(s)^{-1/2}\, ds\, dt.
\end{align*}
Theorem~\ref{thm-decomposition} only guarantees this formula for 
$b$ and $f$ of a particular form. However, functions of these forms
are dense in $C_c(K\times P,A)$ and $C_c(H\times P,A)$ for the 
inductive limit topologies, which are stronger than the topologies 
arising from the imprimitivity bimodule structure, and we can extend the formulas 
to these submodules by continuity.
\end{proof}

\begin{remark}\label{cfaHRW} 
In \cite[Lemma~4.8]{hrw} we used {\it ad hoc} methods to 
find a corresponding isomorphism  for bimodules over the full 
crossed products; there, our formula for $\Omega(\flat(\xi)\otimes b)(p)$ is 
\begin{equation}\label{eq-hrw}
\int_K\int_H\beta_s^{-1}\alpha_t^{-1}
\big( ( f(t)(t\cdot p\cdot s^{-1}))^*b(s)(t\cdot p\cdot s^{-1})\big)
\Delta_K(t)^{-1/2}\Delta_H(s)^{-1/2}\, ds\, dt,
\end{equation}
where $\flat(\xi)$ is the element of the dual module of a crossed product
$W\rtimes K$ based on the $(C\rtimes_\tau H)$--$(\Ind \beta)$ bimodule $W$ 
of \cite{rw85}. To reconcile the formulas, we verify that 
 $\psi:(W\rtimes K){\,}
\widetilde{\;}\to X\rtimes K$ given by 
$\psi(\flat(\xi))(t)=\sigma_t(\xi(t^{-1})^*)\Delta_K(t)^{-1}$ 
is an isomorphism of  imprimitivity bimodules. Then plugging 
$f(t)=\sigma_t(\xi(t^{-1})^*)\Delta_K(t)^{-1}$ into \eqref{eq-must-resolve}
gives \eqref{eq-hrw}, as desired. In retrospect, it seems amazing that in \cite{hrw} we came up 
with the ``right'' formula by the brute-force methods we were using.
\end{remark}
 
\begin{example}\label{ex-green}
It is now easy to see that Definition~\ref{defn-pr} is not symmetric. 
With the notation as above, take $A=\C$, so that $X=\overline{C_c(P)}$ is 
a $C_0(P/H)$--$C_0(P)\rtimes_{\rt}H$-imprimitivity bimodule. 
Lemma~\ref{checkproper} says that the 
action $\sigma=\lt$ of $K$ on $X$ is proper. 
If the corresponding action $\sigma$ of $K$ on $\widetilde X$ 
is proper, then Theorem~\ref{thm-decomposition}(2) 
implies that the action $\sigma$ of $K$ on $C_0(P/H)$ is proper, 
and hence that the action of $K$ on $P/H$ is proper. But 
there are  commuting free and proper actions of $H$ and 
$K$ on $P$ such that $K$ does not act properly on $P/H$: for example, taking
$P=\R$, $H=\Z\subset\R$ and $K=\theta\Z$ for some irrational number in 
$(0,1)$ yields the action of $K\cong\Z$ by powers of an irrational rotation on $\T=\R/\Z$. 
\end{example}

\section{The symmetric imprimitivity theorem for graph algebras}
 
A directed graph $E$ consists of countable sets $E^0$ of vertices and $E^1$ of edges, and range and source maps $r,s:E^1\to E^0$. A Cuntz-Krieger $E$-family in a $C^*$-algebra $A$ consists of partial isometries $\lbrace s_e:e\in E^1\rbrace$ with mutually orthogonal ranges and mutually orthogonal projections $\lbrace p_v:v\in E^0\rbrace$ such that
\[
s_e^*s_e=p_{r(e)},\ \ s_es_e^*\leq p_{s(e)},\text{\ and\ }p_v=\sum_{s(e)=v} s_es_e^*
\text{\ whenever\ } 0<|s^{-1}(v)|<\infty.
\] 
The graph $C^*$-algebra $C^*(E)$ is generated by a universal Cuntz-Krieger family $\{s_e,p_v\}$ (see \cite{KPR} or \cite{RSz}, for example). We write $E^*$ for the path space of $E$, and for $\mu\in E^*$ of length $|\mu|$ we write $s_\mu:=s_{\mu_1}s_{\mu_2}\dots s_{\mu_{|\mu|}}$. The Cuntz-Krieger relations imply that every word in the $s_e$ and $s_f^*$ collapses to one of the form $s_\mu s_\nu^*$ for $\mu,\nu\in E^*$, and these are zero unless $r(\mu)=r(\nu)$. Thus 
\[X_0(E):=\sp\{s_\mu s_\nu^*:\mu,\nu\in E^*, r(\mu)=r(\nu)\}
\]
is a dense $*$-subalgebra of $C^*(E)$.

Suppose we have a left action of a (discrete) group $G$ on $E$ which is free on $E^0$ (and hence is free on $E^1$). The universal property of $C^*(E)$ implies that there is an induced action $\alpha:G\to\Aut C^*(E)$ such that $\alpha_g(s_e)=s_{g\cdot e}$ and $\alpha_g(p_v)=p_{g\cdot v}$. It is shown in \cite[\S1]{PR} that the action $\alpha$ is proper and saturated with respect to $X_0(E)$. Indeed, it is proved in \cite[Lemma~1.1]{PR} that averaging over $\alpha$ gives a linear map $I_G:X_0(E)\to M(X_0(E))^\alpha$ whose range spans the generalized fixed-point algebra $C^*(E)^\alpha$, and that there is an isomorphism $\phi_G$ of the $C^*$-algebra $C^*(G\backslash E)$ of the quotient graph onto $C^*(E)^\alpha$; thus it follows from Rieffel's theory that $C^*(E)\rtimes_{\alpha,r}G$ is Morita equivalent to $C^*(G\backslash E)$. The
maps $I_G$ and $\phi_G$ are also used in \cite{PR} to directly construct a bimodule implementing a symmetric imprimitivity theorem for the full crossed products, as follows. To make cross-referencing easier, we have used the notation of \cite{PR} rather than that of \S\ref{sec-sit}.

Suppose we have commuting free actions of $G$ and $H$ on the left and right of $E$. Because the actions commute, they induce actions on the quotient graphs, and hence we have actions $\alpha:G\to \Aut C^*(E/H)$ and $\beta: H\to \Aut C^*(G\backslash E)$ on their $C^*$-algebras; it is safe to also use $\alpha$ and $\beta$ for the actions on $X_0(E)$ and $C^*(E)$, because the maps $\phi_H$ and $\phi_G$ are then equivariant. For $b\in k(H,X_0(G\backslash E))$, $c\in k(G,X_0(E/H))$ and $x,y\in X_0(E)$, we  define
\begin{align}
b\cdot x&=\sum_{h\in H}\phi_G(b(h))\beta_h(x)\label{graph-la}\\
x\cdot c&=\sum_{g\in G} \alpha_g^{-1}(x\phi_H(c(g)))\label{graph-ra}\\
{}_{k(H,X_0(G\backslash E))}\langle x\,,\,y\rangle(h)&=\phi_G^{-1}\circ I_G(x\beta_h(y^*))\label{graph-lip}\\
\langle x\,,\,y\rangle_{k(G,X_0(E/H))}(g)&=\phi_H^{-1}\circ I_H(x^*\alpha_g(y))\label{graph-rip}
\end{align}
Now \cite[Theorem~2.1]{PR} says that $X_0(E)$ completes to give a Morita equivalence $Z$ between $C^*(G\backslash E)\rtimes_\beta H$ and $C^*(E/H)\rtimes_\alpha G$.

As in the previous section, we aim to apply Theorem~\ref{thm-decomposition} with $X=\overline{X_0(E)}$ the $C^*(G\backslash E)$--$(C^*(E)\rtimes_\alpha G)$ imprimitivity bimodule obtained by ignoring the action of $H$ in \eqref{graph-la}--\eqref{graph-rip}. By \cite[Corollary~3.3]{kqr} or \cite[Corollary~3.1]{PR}, we have \[C^*(E)\rtimes_{\alpha}G=C^*(E)\rtimes_{\alpha}G,\] so we can view $X$ as a module over the reduced crossed product. 

\begin{lemma}
The action $\beta$ of $H$ on $X_0(E)\subset C^*(E)$ induces actions $\beta$ of $H$ on $C^*(G\backslash E)$, $\beta$ on $X$ and $\beta\rtimes\id$ on $C^*(E)\rtimes_\alpha G$, and $(X,H,\beta)$ is then a Morita equivalence between $(C^*(G\backslash E),H,\beta)$ and $(C^*(E)\rtimes_\alpha G,H,\beta\rtimes\id)$. The action $\beta$ on $X$ is proper and saturated with respect to the pre-imprimitivity bimodule ${}_{X_0(G\backslash E)}X_0(E)_{k(G,X_0(E))}$.
\end{lemma}

\begin{proof}
That $\beta$ induces the actions on $C^*(G\backslash E)$ and  $C^*(E)\rtimes_\alpha G$ is standard. Because $\beta$ is compatible with the maps $\phi_G$ and $I_G$ \cite[Lemma~1.7]{PR}, it is easy to check that $\beta$ is compatible with the module actions and inner products. In particular, this implies that each $\beta_h$ is isometric, and hence extends to an action on $X$ implementing the desired Morita equivalence of systems. For the submodule $X_0(E)$, the functions in parts (1) and (2) of Definition~\ref{defn-pr} have finite support, and hence are trivially integrable. For $x,y\in X_0(E)$, the function $\langle x\,,\,y\rangle_D:G\to M(C^*(E))$ defined by
\[
\langle x\,,\,y\rangle_D(g)=I_H(x^*\alpha_g(y))
\]
also has finite support; the embedding  of $M(C^*(E))\times_{\alpha}G$ in $M(C^*(E)\rtimes_\alpha G)$ carries this function into a multiplier $\langle x\,,\,y\rangle_D$ of $k(G,X_0(E))$ which satisfies Definition~\ref{defn-pr}(3). Thus the action of $H$ is proper. To see that it is saturated, we use \cite[Lemma~1.4]{PR} to see that the function $\delta_hs_{G\cdot\mu}s_{G\cdot\nu}^*$ in $k(H,X_0(G\backslash E))$ is given by
\[
\delta_hs_{G\cdot\mu}s_{G\cdot\nu}^*=\delta_h\phi_G^{-1}\circ I_G(s_\mu s_\nu^*)=\langle x\,,\,y\rangle_D
\]
when $x=s_\mu s_\nu^*$ and $y=p_{s(\nu)\cdot h}$.
\end{proof}

Applying Theorem~\ref{thm-decomposition} to $(X,H,\beta)$ gives a $(C^*(E)\rtimes_{\alpha\times\beta}(G\times H))$--$D$ imprimitivity bimodule $\overline{B_1}$, a $(C^*(G\backslash E)\rtimes_\beta H)$--$D$ imprimitivity bimodule $\overline{X_0}$, and a decomposition isomorphism. The space $X_0(E)$ underlies both $\overline{X_0}$ and the bimodule $Z$ of \cite{PR}. Here $X_0$ is really a bimodule over $k(H,X_0(G\backslash E))$ and the generalized fixed-point algebra $D\subset M(C^*(E)\rtimes_\alpha G)$; when we use $\phi_H\rtimes\id$ to identify $D$ with $C^*(E/H)\rtimes_\alpha G$, our formulas convert to the ones \eqref{graph-la}--\eqref{graph-rip} used in \cite{PR}. Thus:

\begin{thm}\label{sitgraphred}
The bimodule ${}_{k(H,X_0(G\backslash E))}X_0(E)_{k(G,X_0(E/H))}$ described in \textnormal{(\ref{graph-la})--(\ref{graph-rip})} completes to give an imprimitivity bimodule which implements a Morita equivalence between $C^*(G\backslash E)\rtimes_{\beta,r} H$ and $C^*(E/H)\rtimes_{\alpha,r} G$.
\end{thm}

Comparing this bimodule to the one for the full crossed products allows us to settle a question left open in \cite[Remark~3.2]{PR}. 

\begin{cor}
Suppose that $G$ and $H$ act freely on the left and right of a directed graph $E$, and let $\alpha$ and $\beta$ denote the induced actions on $C^*(E/H)$ and $C^*(G\backslash E)$. Then regular representations of $(C^*(E/H),G,\alpha)$ are faithful  if and only if regular representations of $(C^*(G\backslash E),H,\beta)$ is faithful.
\end{cor}

\begin{proof}
Let $I$ be the kernel of the quotient map from $C^*(E/H)\rtimes_{\alpha} G$ to $C^*(E/H)\rtimes_{\alpha,r} G$. Then, by the Rieffel correspondence \cite[\S3.3]{tfb}, there are a closed submodule $W$ of the bimodule $Z$ of \cite{PR} and an ideal $J=Z\dashind I$ in $C^*(G\backslash E)\rtimes_{\beta} H$ such that $Z/W$ is a $(C^*(G\backslash E)\rtimes_{\beta} H)/J$--$(C^*(E/H)\rtimes_{\alpha} G)/I$ imprimitivity bimodule. In particular, this implies that the semi-norms on $X_0(E)\subset Z$ induced by the quotient norms on $(C^*(G\backslash E)\rtimes_{\beta} H)/J$ and $(C^*(E/H)\rtimes_{\alpha} G)/I$ coincide \cite[Proposition~3.11]{tfb}. The semi-norm coming from the right inner product is that induced by the reduced norm on $k(G,X_0(E/H))$. However, we know by applying \cite[Proposition~3.11]{tfb} to the bimodule of Theorem~\ref{sitgraphred} that this coincides with the seminorm induced by the left inner product and the reduced norm on $k(H,X_0(G\backslash E))$. Thus the seminorm on $k(H,X_0(G\backslash E))$ pulled back from the quotient $(C^*(G\backslash E)\rtimes_{\beta} H)/J$ is the reduced seminorm, the quotient is the reduced crossed product, and $J$ is the kernel of the quotient map onto $C^*(G\backslash E)\rtimes_{\beta,r} H)$. Since $I=0$ if and only if $Z\dashind I=0$, the result follows.
\end{proof}

To complete the analysis, we identify $\overline{B_1}$ and check that the decomposition 
of Theorem~\ref{thm-decomposition} gives an isomorphism between the tensor-product equivalence of \cite[Theorem~1.4]{PR} 
and that of Theorem~\ref{sitgraphred}. 

The algebra $B_1$ is spanned by the range of the inner product on $X_0(E)$, and since 
\[
\delta_g s_\mu s_\nu^*=\langle s_\nu s_\mu^*\,,\,p_{g^{-1}\cdot s(\nu)}\rangle_{k(H,X_0(E))},
\]
this is all of $k(H,X_0(E))$. This is also a dense subspace of the Combes bimodule $Y\rtimes_\alpha G$, where $Y$ is the $(C^*(E)\rtimes_\beta H)$--$C^*E/H)$ imprimitivity bimodule obtained by ignoring $G$ in \eqref{graph-la}--\eqref{graph-rip}. To see that $\overline{B_1}$ is isomorphic to $Y\rtimes_{\alpha,r} G$, we need to note that the map $\phi_H\rtimes \id$ is an isomorphism of $C^*(E/H)\rtimes_{\alpha,r} G$ onto the fixed-point algebra $D$, and check that the inner products and module actions match up. This follows from calculations like those done in the proof of Theorem~\ref{astridsiso}. The isomorphism is given by
\[
\Omega(f\otimes e)=\sum_{h\in H} f(h)\cdot (\beta\rtimes\id)_h(e),
\]
where the action is that of $E_0\subset k(G,X_0(E/H))$ on $X_0$. Working out the formulas in terms of the product in $X_0(E)\subset C^*(E)$ gives
\begin{equation}\label{decompforgraphs}
\Omega(f\otimes e)=\sum_{h\in H}\sum_{g\in G} \alpha_g^{-1}\big(f(h)\phi_H(\beta_h(e(g)))\big)
\end{equation}
This time the functions $f,e$ of the required form span $k(H,X_0)$ and $k(G,Y_0)$, respectively. Hence:

\begin{cor}
The map $\Omega:k(H,X_0)\otimes k(G,Y_0)\to X_0$ defined by \textnormal{(\ref{decompforgraphs})} extends to an isomorphism of $(X\rtimes_{\beta,r}H)\otimes_{C^*(E)\rtimes_r(H\times G)}(Y\rtimes_{\alpha,r}G)$ onto the reduced symmetric imprimitivity bimodule of Theorem~\ref{sitgraphred}.
\end{cor}

\end{document}